%% file: KPK-MMS2011.tex
\documentclass[a4paper,11pt,oneside,abstracton,DIV14,final]{scrartcl}

\usepackage[english]{babel}
\usepackage[latin1]{inputenc}
\usepackage[T1]{fontenc}

\usepackage{mstyle}

\usepackage{graphicx}
\usepackage{subfigure}
\usepackage{booktabs}

\usepackage{color}

\usepackage[ruled]{algorithm} 
\usepackage{algorithmic}

\title{Semi-Parametric Drift and Diffusion Estimation for Multiscale Diffusions}

\author{S.~Krumscheid\footnotemark[2]
\and G.~A.~Pavliotis\footnotemark[2]\ \footnotemark[4]\ \footnotemark[5]
\and S.~Kalliadasis\footnotemark[3]}
\date{February $19$, $2013$}

\begin{document}
\maketitle

\renewcommand{\thefootnote}{\fnsymbol{footnote}}

\footnotetext[2]{Department of Mathematics, Imperial College London, South Kensington, London SW7 2AZ, United Kingdom: \{s.krumscheid10, g.pavliotis\}@imperial.ac.uk}
\footnotetext[3]{Department of Department of Chemical Engineering, Imperial College London, South Kensington, London SW7 2AZ, United Kingdom: s.kalliadasis@imperial.ac.uk}
\footnotetext[4]{INRIA MicMac Project Team, Ecole des Ponts ParisTech, 6 et 8 Avenue Blaise Pascal, Cit{\'e} Descartes - Champs sur Marne, 77455 Marne la Valle Cedex 2, France: grigorios.pavliotis@cermics.enpc.fr}
\footnotetext[5]{Institut f{\"u}r Mathematik, Freie Universit{\"a}t Berlin, Arnimallee 6, 14195 Berlin, Germany: gpavlo@zedat.fu-berlin.de}

\renewcommand{\thefootnote}{\arabic{footnote}}

\begin{abstract}
  We consider the problem of statistical inference for the effective dynamics of multiscale
  diffusion processes with (at least) two widely separated characteristic time scales. 
  More precisely,
  we seek to determine parameters in the effective equation describing
  the dynamics on the longer diffusive time scale, i.e.\ in a homogenization framework.
  We examine the case where both the drift and the diffusion coefficients in the effective dynamics are
  space-dependent and depend on multiple unknown parameters.
  It is known that classical estimators, such as Maximum Likelihood and Quadratic
  Variation of the Path Estimators, fail to obtain reasonable estimates for parameters in the effective
  dynamics when based on observations of the underlying multiscale diffusion. We propose a
  novel algorithm for estimating both the drift and diffusion coefficients in the effective dynamics based
  on a semi-parametric framework.  We demonstrate by means of extensive numerical simulations of a number of
  selected examples that the algorithm performs well when applied to data from a multiscale diffusion. These
  examples also illustrate that the algorithm can be used effectively to obtain accurate and unbiased estimates.
\end{abstract}



\input{introduction}

\input{estimator}

\input{numerics}

\input{outlook}

\section*{Acknowledgements}
We are grateful to the anonymous referees for insightful comments
and suggestions. We also thank R.\ N{\"u}rnberg and A.\ M.\ Stuart
for critically reading an earlier version of the manuscript and for
their useful comments and suggestions. This work is supported by the
Engineering and Physical Sciences Research Council of the UK through
Grant No.\ EP/H034587.

\bibliography{%
  ./bibs/SDE,%
  ./bibs/NumSimSDE,%
  ./bibs/EstimatesInSDEs,%
  ./bibs/generalMath,%
  ./bibs/Homogenization,%
  ./bibs/PhysicsSDEs,%
  ./bibs/SPDE}

\bibliographystyle{siam}

\end{document}

%% file: introduction.tex
\section{Introduction}
\label{sec:intro} Problems with multiple temporal and/or spatial
scales emerge naturally in a wide variety of fields in science and
engineering. From biological phenomena \cite{Chauviere2010} and
atmosphere and ocean science \cite{Majda2008} to molecular dynamics
\cite{Griebel2007}, material science \cite{Fish2009}, fluid and
solid mechanics \cite{Horstemeyer2010,Huerre1998}. Many of such
systems are often subject to noise that is either due to thermal
fluctuations \cite{Einstein1956}, randomness in the environment
(e.g.\ uncertainty in some parameters)
\cite{Horsthemke1984,Savva2010,Pradas2011}, coarse-graining of
high-dimensional deterministic systems with random initial
conditions \cite{Mazo2002,Zwanzig2001}, or stochastic
parameterization of small scale effects \cite{E2005}.

Mathematically the influence of noise in a system can be described
by a single or coupled, nonlinear stochastic differential equations
(SDEs) with multiple scales. In many cases these equations are
high-dimensional but from a practical point of view only the
evolution of some components of the solution is of main interest,
since they act on a slower scale. It is therefore desirable to
approximate the full system by an adequate simplified
low-dimensional effective model (coarse-graining) that retains the
essential dynamic characteristics of the full system. The effective
equation is often amenable to analytical and numerical work.
However, usually only the complete multiscale (fast/slow) system is
directly observable but not the effective dynamics. Consequently,
much research has been undertaken to find accurate approximations to
the effective dynamics \cite{Givon2004}.

In many applications a wealth of data (i.e.\ observations of the
fast/slow system) is often available and it is therefore worthwhile
to use these data to determine the effective dynamics. However,
extracting the effective dynamics directly from the available data
is often not straightforward and it is thus important to adapt
techniques from analysis, statistics, and numerical analysis to
obtain appropriate coarse-grained models.

This data-driven coarse graining methodology has been applied to
relatively simple (stochastic) systems for which a low dimensional
coarse-grained equation exists, using techniques such as
homogenization and averaging. For example, in \cite{Pavliotis2007}
Brownian motion in a two-scale potential was studied. Therein,
through both rigorous mathematical analysis and numerical
simulations, it was shown that the estimation of drift and diffusion
coefficients in the coarse-grained model is asymptotically biased
when using classical estimators. Furthermore, it was shown that
subsampling the available data at an appropriate rate between the
two characteristic time scales of the full system is necessary for
an accurate estimation of both drift and diffusion coefficients in
the coarse-grained model. More general fast/slow systems of SDEs for
which a coarse-grained equation exists were studied in
\cite{Papavasiliou2009} where it was shown that the same issue of
asymptotically biased estimators persists in the homogenization
framework -- that is, when considering an effective dynamics on the
(longer) diffusive time scale -- and that appropriately subsampled
data reduce the bias. These techniques were then applied to the
problem of estimating eddy diffusivities from noisy Lagrangian
observations in \cite{Cotter2009} where an improved algorithm that
combines subsampling with appropriate averaging and variance
reduction techniques was proposed and tested. Furthermore, inverse
problems for multiscale partial differential equations -- a problem
closely related to that of parameter estimation -- were studied in
\cite{Nolen2012}.

Although appropriate subsampling can potentially improve the
accuracy of estimators in the context of multiscale diffusions, as it
yields unbiased estimators, the
question of an optimal subsampling rate (i.e.\ the rate for which
the biased is removed) remains open; e.g. the studies in
\cite{Pavliotis2007,Papavasiliou2009} provide only existence
results. Furthermore, the numerical experiments presented in the
aforementioned works indicate that the optimal subsampling rate is
not only problem dependent, but is also different for different
parameters in the same model. Consequently, an optimal subsampling
rate is in general unknown. The problem of identifying an optimal
subsampling rate for Gaussian processes has been studied in
\cite{Azencott2010,Azencott2011}; see also
\cite{Zhang2005} for related work in the context of econometrics.\\

Related problems have been studied in the context of numerical
analysis for SDEs with multiple time scales. In particular, the
heterogeneous multiscale method (HMM)
\cite{Vanden-Eijnden2003,E2005} is based on the idea of evolving the
solution of the low-dimensional coarse-grained equation, when the
coefficients in the coarse-grained equation are being evaluated ``on the
fly'' by running short runs of the underlying fast dynamics. Similar
ideas have been proposed in the framework of the ``equation-free
methodology'' introduced by Kevrekidis and collaborators; see e.g.\
\cite{Theodoropoulos2000,Kevrekidis2003,Kevrekidis2004,Kevrekidis2009}.
Recently, the HMM methodology has been extended to approximate
stochastic partial differential equations with multiple timescales
\cite{Abdulle2012}. As such, these techniques can be considered as a
hybrid between numerical analysis and statistical inference since
the coefficients in the coarse grained SDE are estimated from data
that are obtained from short runs of the full dynamics. We
emphasize, however, that in the HMM the fast dynamics is assumed to
be known, whereas in the aforementioned works on parameter
estimation for multiscale diffusions, as well as in the present work
no such assumption is made.

We also mention that the effect of the multiscale structure on the evolution of the
coarse-grained probability density using the Fokker-Planck equation
(Kolmogorov's forward equation) was studied in \cite{Frederix2011}. In
this study it was shown that when decreasing the spatial
discretization in a finite difference approximation the error
increases rapidly and that in order to avoid this, it is necessary
to improve the accuracy of the estimators of the drift and diffusion
coefficients.

In many cases of interest the noise in the coarse-grained equation appears in a multiplicative way. 
One is thus confronted with the problem of estimating parameters in both the 
drift and the diffusion coefficients of an SDE of the form
\begin{align*}
dx_t = f(x_t;\vartheta)\,dt + g(x_t;\theta)\,dW_t\;,
\end{align*}
with unknown parameters $(\vartheta,\theta)^T\in\Theta$, the set of
all feasible parameters. The problem is further complicated by the
fact that the parameters $\vartheta$ and $\theta$ may not be
independent as for example in the problem of Brownian motion in a
two-scale potential, see Section~\ref{sec:num:multi:langevin}.
In the absence of multiscale effects in the data, or if one
assumes that the optimal subsampling rate is known, a combination of
the maximum likelihood estimator (MLE) and the quadratic variation
of the path (QVP) is the most commonly used estimator in practice;
see \cite{PrakasaRao1999,Kutoyants2004,Liptser2010} for background material
on classical estimators. The QVP (or a variant thereof) is used to
estimate parameters in the diffusion coefficient and, based on these
estimates, the MLE is used to obtain estimates for the parameters in
the drift. As the MLE is based on an SDE with unit diffusion, one
usually transforms the original SDE into an SDE with unit diffusion
coefficient by applying It{\^o}'s formula to
\begin{align}
  h_\theta(x) = \int_c^x g(u;\theta)^{-1}\,du\;,\label{eq:trans:lamperti}
\end{align}
for an arbitrary $c$ in the state space of the process $x$ and
replacing $\theta$ with the estimator obtained from the QVP for
concreteness. For example, when considering
\begin{align*}
  g(x;\theta) = \sqrt{\theta_1 + \theta_2x^2}\;,
\end{align*}
which is the diffusion coefficient in the coarse-grained equation of the 
stochastic Burgers' equation (see Section \ref{sec:num:multi:burgers}),
then transformation \eqref{eq:trans:lamperti} reads
\begin{align*}
  h_\theta(x) = \frac{\ln{\bigl(\sqrt{\theta_2}x+\sqrt{\theta_1+\theta_2x^2}\bigr)} -
  \ln{\bigl(\sqrt{\theta_2}c+\sqrt{\theta_1+\theta_2c^2}\bigr)}}{\sqrt{\theta_2}}\;,
\end{align*}
where $\theta=(\theta_1,\theta_2)^T$ has to be replaced by a
previously obtained estimator. 

Notice that this transformation can be singular when implementing it in practice,
so that special care has to be taken when performing numerical simulations.
In particular within the MLE framework where a (nonlinear) objective function
needs to be maximized, this might cause problems. An
alternative approach to estimate parameters is based on the MLE for
the discretized approximation, e.g.\ obtained by the Euler-Maruyama
scheme \cite[ch.\ $9.1$]{Kloeden1992}, of the time-continuous SDE.
However, the non-vanishing (fixed) time-step size introduces an
additional bias so that even for simple models this approach does
not necessarily yield consistent estimators \cite{Lo1988}.\\

The main aim of the present study is to develop statistical
inference techniques that enable us to estimate parameters in
both drift and diffusion coefficients of a coarse-grained equation in the
presence of an underlying (either stochastic or deterministic)
multiscale structure in the fast/slow system. 
More precisely, given only observations of the slow component
of the fast/slow system without any further knowledge of the fast
component, the aim is to infer the coefficients in the coarse-grained equation.
Furthermore, we wish
to extend this approach in a semi-parametric framework, to
situations where the drift and diffusion coefficients can be
expanded in an appropriate (e.g.\ Taylor series) expansion.

Besides the lack of reliable statistical inference techniques for
these problems, the motivation for this study originates also from
recent results on the derivation of coarse-grained equations (also
known as amplitude equations) for stochastic partial differential
equations (SPDEs) with quadratic nonlinearities \cite{Blomker2007}.
A typical example for such an SPDE is the stochastic Burgers'
equation
\begin{align}
du_t = \Bigl((\partial_x^2 + 1)u_t + \partial_xu_t^2 + \varepsilon^2u_t \Bigr)\,dt
  + \varepsilon \mathcal{Q}\,d\mathcal{W}_t\;,\label{eq:spde:generic}
\end{align}
on $[0,\pi]$ equipped with appropriate boundary conditions.
Therein $\mathcal{Q}$ denotes the covariance operator, $\mathcal{W}$
space-time white noise, and $0<\varepsilon \ll 1$. 
To study solutions to \eqref{eq:spde:generic} of $O(\varepsilon)$ on time 
  scales of $O(\varepsilon^{-2})$, i.e.\ ensuring that we 
  are in the regime described by amplitude equations, a diffusive rescaling is
  performed by defining $v$ via $\varepsilon v(\varepsilon^2t)=u(t)$.
Then $v$ solves
\begin{align}
  dv_t = \Bigl(\frac{1}{\varepsilon^2}(\partial_x^2 + 1)v_t + \frac{1}{2\varepsilon}\partial_xv_t^2 + v_t \Bigr)\,dt
  + \frac{1}{\varepsilon} \mathcal{Q}\,d\mathcal{W}_t\;.
  \label{eq:spde:generic:diffusive}
\end{align}
If the SPDE is equipped with homogeneous Dirichlet boundary conditions it can be shown that the dominant
mode of the solution to \eqref{eq:spde:generic:diffusive} can be approximated by the solution to an one-dimensional
SDE driven by classical Brownian motion $W_t$ of the form
\begin{align*}
  dX_t = \bigl(AX_t-BX_t^3\bigr)\,dt + \sqrt{\sigma_a^2 + \sigma_b^2X_t^2}\,dW_t\;.
\end{align*}
Since this class of SPDEs arises in many different applications --
from population biology~\cite{Horsthemke1984} to fluid
dynamics~\cite{Pradas2011,Pradas2012} where many data are available
-- it is of major interest to obtain effective dynamics for the
dominant modes of solutions to the SPDE by means of a data-driven
coarse graining methodology.

The statistical inference technique we propose here consists of two
steps. First, we use the martingale property of the stochastic
integral to obtain an equation involving only the drift but not the
diffusion coefficient of the SDE. Since the drift might depend
on multiple (unknown) parameters, it is generally impossible to
obtain the parameters uniquely from a single equation. The
main element we employ to overcome this under-determined situation
is the often disregarded initial condition. In fact, by varying the
initial condition of the SDE one can define the estimator for the
drift parameters via the best approximation of a system of
equations. For the second step, we rely on the estimators for the
drift parameters and on the It{\^o} isometry to obtain a relation
for the unknown parameters in the diffusion. Using the same idea as
in the first step, that is by varying the initial condition of the
SDE, we can also define the estimators for parameters concerning the
diffusion via the best approximation of a system of equations. 
The expectations involved in both steps of
the estimation procedure are approximated by an average over many
short trajectories (i.e.\ by ensemble averages) in contrast to 
classical estimators that often rely on a long trajectory 
of the underlying process, see e.g.\ \cite{PrakasaRao1999}.
That is, the methodology we propose here relies on independent short
  trajectories starting from different initial conditions.\\

The main practical advantages of the methodology proposed here can
be summarized as follows:
\begin{itemize}
\item[(a)] 
Coarse-grained models of high-dimensional, possibly infinite-dimensional,
problems that retain the essential dynamic characteristic of the full
system are a powerful tool to study problems arising in science and
engineering. The simplicity of such an effective model, typically an SDE,
makes it attractive for mathematical and numerical scrutiny.
When the coefficients of the SDE can be computed exactly, the effective
SDE allows us to decipher rapidly some of the basic characteristics of 
the full system, e.g.\ by computing first-passage properties and properties alike
as it was done for instance in \cite{Pradas2011,Pradas2012}.
However, in several problems it is not straightforward to obtain
the coefficients of the coarse-grained SDE exactly due to the complexity of
the underlying full system. From a theoretical standpoint several 
assumptions need to be made, e.g.\ the dynamics is dominated by a single 
eigenfunction/mode (often associated with a symmetry in the system) 
and the higher-order modes decay sufficiently fast (see the derivation of
the ``phase-diffusion'' equation describing the transverse
instability of propagating waves/fronts; e.g.\ in
\cite{Kuramoto2003,Kalliadasis2000}), to obtain computable coefficients.
Conversely, from a practical point of view obtaining approximations (estimators) 
of the coefficients in the coarse-grained model, which can be used
to study the basic characteristics of the full system, is very appealing.
\item[(b)] A model for a physical or technological process might not be readily
available, either because its derivation is cumbersome, or it is not
straightforward to formulate it from first principles. But the
underlying physics of the phenomena at hand and previous experience
with similar systems suggests an SDE of the form adopted here, at
least in certain regions of the parameter space. In a spirit similar
to the equation-free approach, one can utilize available data and
obtain a low-dimensional approximation to the process, which in turn
can be used as a model for the process in regions of the parameter
space consistent with the assumptions imposed from the outset. But
often such models can be applicable in regions beyond those dictated
by the assumptions. For instance, for a long-wave instability such
as that observed on a surface of a film flowing down an inclined
plane \cite{Kalliadasis2012} the growth rate curve of infinitesimal
disturbances extends from the origin up to a cut-off wavenumber. The
Landau-Stuart equation is then only applicable sufficiently close to
the cut-off, i.e.\ it can be used to describe the transition when an
unstable wave motion of given (``fundamental'') wavenumber interacts
with its first stable harmonics. On the other hand, when the growth
rate curve has a ``nose'' near criticality consistent with a
short-wave instability and so that the system equilibrates to a
stationary norm solution in the nonlinear regime, such as with with
Rayleigh-B\'enard convection \cite{Cross1993}, the Landau-Stuart
equation can be applicable even past the instability threshold.
\end{itemize}

The rest of the paper is organized as follows. In Section
\ref{sec:estimator} we present the precise derivation of the
estimators for systems with and without multiscale structure. Based
on these results, in Section \ref{sec:num} the general applicability
and performance of the prosed methodology is investigated via
different numerical examples. The extensive numerical study we
undertake illustrates that the proposed technique enables us to
estimate accurately parameters in multiscale diffusions. A summary
of the results and future perspectives are offered in Section
\ref{sec:outlook}.


%% file: estimator.tex
\section{Estimators}
\label{sec:estimator}

We present the precise derivation of the drift and diffusion
estimators for systems without and systems with multiscale effects
present. First we outline the methodology for SDEs without
multiscale structure and illustrate some properties of the
estimators in this case, before presenting the set-up for multiscale
diffusions.

\subsection{Derivation of Drift and Diffusion Estimators}
\label{sec:estimator:derivation}

For the sake of simplicity we consider here only one dimensional
real-valued processes (see Sections \ref{sec:num:multi:langevin2d} and
\ref{sec:num:multi:taylorgreen} for examples of multivariate processes).
Consider
the scalar-valued It{\^o} stochastic differential
equation
\begin{align}
  dx_t = f(x_t)\,dt + \sqrt{g(x_t)}\,dW_t\;,\quad x(0) \equiv x_0 = \xi\;,\label{eq:ampl:sde}
\end{align}
where $W$ denotes standard one-dimensional Brownian motion. Both the drift coefficient
$f$ and the diffusion coefficient $g$
are assumed to be sufficiently smooth such that the SDE provides a
unique solution for any initial condition $\xi\in\R$ and on
any finite time interval; see e.g.\ \cite[sec.\ $1$]{Krylov1999}.
Moreover, we assume that both drift $f$ and diffusion $g$ depend on
unknown parameters and the task is to estimate the parameters in $f$
and $g$ from available data. In fact, here we focus on the case when
$f(x)$ and $g(x)$ are polynomials in $x$ of degree $\max{\{J_f\}}$ and
$\max{\{J_g\}}$, respectively, where $J_f,J_g \subset \N_0$ denote
index sets of finite cardinality $p=\abs{J_f}$ and $q=\abs{J_g}$
respectively. The unknown coefficients of the polynomials are
$\vartheta\equiv{(\vartheta_j)}_{j\in J_f}\in\R^{p}$ and
$\theta\equiv {(\theta_j)}_{j\in J_g}\in\R^{q}$ respectively, i.e.\
we consider
\begin{align}
f(x)\equiv f(x;\vartheta) := \sum_{j\in J_f} \vartheta_j x^j\quad\text{and}\quad
g(x)\equiv g(x;\theta) := \sum_{j\in J_g} \theta_j x^j\;.\label{eq:parameterization}
\end{align}
Consequently, $f$ and $g$ are linear functions in $\vartheta$ and
$\theta$, respectively. These particular assumptions on the drift
$f$ and the diffusion $g$ simplify the notation in what follows and
they will lead to a linear system of equations for the estimators of
the parameters. Notice that parameterization \eqref{eq:parameterization} is not compulsory and also other 
parameterizations, such as $(x;\theta)\mapsto\sum_{j\in J}\theta_jv_j(x)$ for some known functions $v_j$,
will lead to linear system of equations for $\theta$. More general parameterizations will, however, not
necessarily result in a system of linear equations for $\theta$ anymore and more elaborated methods 
(e.g.\ an iterative approach using Newton's method) are required.

The starting point in the derivation of the estimators is based on the following identities
\begin{subequations}
  \begin{align}
    &\E(x_t - \xi) = \int_0^t\E\bigl(f(x_s)\bigr)\, ds\;,\label{eq:est:start:drift}\\
    &\E\Bigl(\bigl(x_t - \xi - \int_0^tf(x_s)\, ds\bigr)^2\Bigr) =
    \int_0^t\E\bigl(g(x_s)\bigr)\, ds\;,\label{eq:est:start:diffusion}
  \end{align}
  \label{eq:est:start}%
\end{subequations}
owing to the martingale property of the stochastic integral and the
It{\^o} isometry, respectively, holding for any fixed initial
condition $\xi$. The next step is to incorporate the
parameterization of the functions $f$ and $g$ into
\eqref{eq:est:start}, to identify the functional structure of the
relation for the parameters $\vartheta$ and $\theta$, respectively.
We begin with substituting the parameterization of $f$ into equation
\eqref{eq:est:start:drift}, which will yield an estimator for the
parameter $\vartheta$, that is present in the drift term alone.
Based on this estimator it is possible to proceed similarly with
equation \eqref{eq:est:start:diffusion} and
eventually obtain an estimator for the parameter $\theta$ present in the diffusion part.\\

Substituting ansatz \eqref{eq:parameterization} for $f$ into \eqref{eq:est:start:drift} yields
\begin{align}
  \E(x_t - \xi) = \sum_{j\in J_f}\vartheta_j \int_0^t \E({x_s}^j)\,ds\;,\label{eq:est:start:drift:param}
\end{align}
for a given initial condition $\xi$. Therein $\E$ denotes the
expectation with respect to the Wiener measure and with respect to
processes starting at a fixed initial condition $\xi$. To emphasize
this dependency on the initial condition we will use the notation
$\E\equiv\E_{\xi}$. Fix a time $t>0$ (the
question how to chose the final time $t$ will be addressed in
Section \ref{sec:num}) and define
\begin{align*}
  b_1\colon\R\ni\xi\mapsto b_1(\xi) &:= \E_\xi(x_t - \xi)\in\R\\
  a_1\colon\R\ni\xi\mapsto a_1(\xi)
  &:= \Bigl(\int_0^t \E_\xi({x_s}^j)\,ds\Bigr)_{j\in J_f}\in\R^p\;.
\end{align*}
With these definitions equation \eqref{eq:est:start:drift} can be rewritten as
\begin{align}
  a_1(\xi)^T\vartheta = b_1(\xi)\;.\label{eq:fun:form:param:single}
\end{align}
The above equation is under-determined for $p>1$. To derive a
well-defined estimator for $\vartheta$, we consider a finite
sequence of initial conditions ${(\xi_{i})}_{1\le i\le m}$ with
$m\ge p$. Since \eqref{eq:fun:form:param:single} is valid for each
initial condition, this approach yields a system of linear equations
\begin{align}
  A_1\vartheta = b_1\;,\label{eq:fun:form:param:sys}
\end{align}
with $A_1 := \bigl(a_1(\xi_{i})^T\bigr)_{1\le i\le m}\in\R^{m\times p}$ and
$b_1:= \bigl(b_1(\xi_{i})\bigr)_{1\le i\le m}\in\R^{m}$.
The linear system does not have a unique solution in general (if a solution exists at all).
To overcome this shortcoming we define the solution of the system of linear equations in
\eqref{eq:fun:form:param:sys}, i.e.\ the estimator of the drift parameter, to be the \emph{best
approximation}
\begin{align}
  \hat{\vartheta} &:= \argmin_{s\in\mathcal{S}_1}\norm{s}_2^2\;,
  \quad\mathcal{S}_1:=\bigl\{z\in\R^p\colon \norm{A_1z-b_1}_2^2 \rightarrow\min\bigr\},\nonumber
\end{align}
respectively,
\begin{align}
  \hat{\vartheta} &:= A_1^{+}b_1\;,\label{eq:drift:est}
\end{align}
with $A_1^{+}$ being the Moore-Penrose pseudo-inverse
\cite{Ben-Israel2003}. We note that the estimation of parameters in
the drift does not require knowledge of the diffusion coefficient.

Assume now that we have already estimated the parameters in the drift $f$ correctly.
Then substituting the ansatz \eqref{eq:parameterization}
of $g$ into \eqref{eq:est:start:diffusion} yields
\begin{align*}
  \E\Bigl(\bigl(x_t - \xi -\int_0^tf(x_s;\hat\vartheta)\, ds\bigr)^2\Bigr) =
  \sum_{j\in J_g}\theta_j \int_0^t \E({x_s}^j)\,ds\;.
\end{align*}
Recall again that the expectation is with respect to the Wiener
measure and processes starting at $\xi$, hence $\E\equiv\E_{\xi}$.
To cope with multiple parameters, we follow the same approach as for the drift
parameters above. In fact, define here
\begin{align*}
  b_2\colon\R\ni\xi\mapsto b_2(\xi) &:=
  \E_\xi\Bigl(\bigl(x_t - \xi-\int_0^tf(x_s;\hat\vartheta)\, ds\bigr)^2\Bigr)\in\R\\
  a_2\colon\R\ni\xi\mapsto a_2(\xi)
  &:= \Bigl(\int_0^t \E_\xi({x_s}^j)\,ds\Bigr)_{j\in J_g}\in\R^q
\end{align*}
and consider again a finite sequence of initial conditions ${(\xi_{i})}_{1\le i\le m}$.
Then we also obtain
a system of linear equations for the parameters
\begin{align}
  A_2\vartheta = b_2\;,\label{eq:fun:form:param:sys2}
\end{align}
with $A_2 := \bigl(a_2(\xi_{i})^T\bigr)_{1\le i\le m}\in\R^{m\times q}$ and $b_2:= \bigl(b_2(\xi_{i})\bigr)_{1\le i\le m}\in\R^{m}$. We define the estimator again via the best approximation
\begin{align}
  \hat{\theta} := A_2^{+}b_2\;.\label{eq:diffusion:est}
\end{align}
Since the estimation of diffusion parameters $\theta$ is based on
the estimators $\hat{\vartheta}$ for the drift parameters,
additional error sources might affect the estimator $\hat{\theta}$
in practice
-- see Section \ref{sec:estimator:properties} for an example of
this error propagation.

In practice we are confronted with discrete observations instead of
continuous ones, so that we need to approximate the (deterministic)
integrals in $a_1(\cdot)$, $a_2(\cdot)$, and $b_2(\cdot)$. Assume
that we have $(n+1)$ observations at equidistant times $t_k := kh$
for $0\le k\le n$ and $h:=t/n$. The goal is to approximate the
integrals by means of these observations. Since the integrands 
depend on the path of the solution of an SDE, we cannot expect the integrands
to be very smooth. For such ``rough'' functions the trapezoidal rule
is more accurate than Simpson's rule \cite{Cruz-Uribe2002}.
Consequently, we approximate the various integrals via the composite
trapezoidal rule
\begin{align}
  \int_0^t \E({x_s}^j)\,ds = \sum_{k=0}^{n-1}\int_{t_k}^{t_k+h}\E({x_s}^j)\,ds
  \approx \frac{h}{2}\Bigl(\E({x_0}^j\bigr) + \E({x_t}^j\bigr) +
  2\sum_{k=1}^{n-1}\E({x_{t_k}}^j)\Bigr)\;,\label{eq:intapprox:trapez}
\end{align}
where we used that $t_0 = 0$ and $t_n = t$.

\subsection{Description of the Algorithm: An Example}
\label{sec:estimator:algorithm}

To apply the actual parametric estimation procedure introduced in the
previous section to a specific problem, not only data need to be
available but also a parameterization needs to be chosen.
Consequently, the complete algorithm can be understood as consisting of two
stages:
\begin{enumerate}
\item \textbf{Initialization:} The time step $h$ is given by the
  underlying time series of observations and is assumed to be
  constant (however, this assumption is not necessary, in
  fact, the procedure might be carried out in the exact same
  manner with a non-equidistant sampling rate) and
  the terminal time $t=nh$ is fixed by choosing $n$ appropriately 
  (cf.\ Section \ref{sec:num}).
  Expectations are approximated by averages over $N$ trajectories generated
  form $N$ independent Brownian motions.
  The crucial step is to fix a
  parameterization for both drift and diffusion coefficients. Lastly, the sequence of initial conditions
  ${(\xi_{i})}_{1 \le i \le m}$ needs to be chosen appropriately.
\item \textbf{Two-step Estimation:} Based on the initializations in the previous stage the estimators are
  well-defined. According to Section \ref{sec:estimator:derivation}
  the parameters $\vartheta$ and $\theta$ are estimated
  successively, first the parameters in the drift and then the parameters in diffusion coefficient.
  For both estimators, two steps need to be performed:
  \begin{enumerate}
  \item \textbf{Assembling} the linear system equations \eqref{eq:fun:form:param:sys} and
    \eqref{eq:fun:form:param:sys2} respectively.
  \item \textbf{Solving} the arising systems via best approximation \eqref{eq:drift:est} and
    \eqref{eq:diffusion:est}
    respectively.
  \end{enumerate}
\end{enumerate}
Since the estimation step depends on the considered parameterizations for drift and diffusion, we present
a detailed pseudocode of the methodology in Algorithm \ref{alg:est}
\begin{algorithm}
\begin{algorithmic}[1]
\REQUIRE $h>0$ and $X\in\R^{m\times N\times (n+1)}$
\FOR{$i=1$ to $m$}
\FOR{$j=1$ to $n$}
\STATE $\alpha_j \leftarrow \sum_{k=1}^N \frac{X_{i,k,j+1}}{N}$
\STATE $\beta_j \leftarrow \sum_{k=1}^N \frac{(X_{i,k,j+1})^3}{N}$
\ENDFOR
\STATE $A_{i,1}\leftarrow\frac{h}{2}(X_{i,1,1}+\alpha_n+2\sum_{j=1}^{n-1}\alpha_j)$
\STATE $A_{i,2}\leftarrow\frac{h}{2}\bigl((X_{i,1,1})^3+\beta_n+2\sum_{j=1}^{n-1}\beta_j\bigr)$
\STATE $b_{i}\leftarrow \alpha_n - X_{i,1,1}$
\ENDFOR
\STATE $\vartheta \leftarrow A^{+}b$
\FOR{$i=1$ to $m$}
\FOR{$j=1$ to $n$}
\STATE $\gamma_j \leftarrow \sum_{k=1}^N \frac{(X_{i,k,j+1})^2}{N}$
\ENDFOR
\STATE $A_{i,1}\leftarrow nh$
\STATE $A_{i,2}\leftarrow\frac{h}{2}\bigl((X_{i,1,1})^2+\gamma_n+2\sum_{j=1}^{n-1}\gamma_j\bigr)$
\FOR{$j=1$ to $N$}
\STATE $\delta_j \leftarrow \frac{h}{2}(f(X_{i,1,1};\vartheta) + f(X_{i,j,n+1};\vartheta) + 2\sum_{k=2}^{n}f(X_{i,j,k};\vartheta))$
\ENDFOR
\STATE $b_{i}\leftarrow \sum_{j=1}^N \frac{(X_{i,j,n+1} - X_{i,1,1} - \delta_j)^2}{N}$
\ENDFOR
\STATE $\theta \leftarrow A^{+}b$
\RETURN $(\vartheta^T,\theta^T)^T$
\end{algorithmic}
\caption{Algorithm for the estimation of the parameters in the drift and diffusion coefficients in
  \eqref{eq:algo:ex:GL}.}
\label{alg:est}
\end{algorithm}
for the example using $J_f = \{1,3\}$ and $J_g = \{0,2\}$, i.e.\ $p=2=q$ with
\begin{align}
  f(x;\vartheta) = \vartheta_1x + \vartheta_3x^3\quad\text{and}\quad g(x;\theta) = \theta_0 + \theta_2x^2\;.
  \label{eq:algo:ex:GL}
\end{align}
This setting corresponds to the Landau-Stuart equation that will
play a vital role in the numerical examples discussed in Section
\ref{sec:num}. The input arguments of Algorithm \ref{alg:est} are the
time step size $h$ and the data array $X$. The
dimension of the array is a result of $m$ different initial
conditions each with $N$ trajectories of $(n+1)$ observations.
When we denote by $x_{t\vert\xi}^{i}$ the value at time $t$ of the
$i$-th trajectory started initially in $\xi$, then $X$ corresponds 
to the collection of these trajectories at discrete (equidistant) times. For the example
we consider here, we define approximations of the first three moments at time $t$ via
\begin{align*}
  \bar{x}_{t\vert\xi} := N^{-1}\sum_{i=1}^{N}x_{t\vert\xi}^{i}\;,\quad
  \tilde{x}_{t\vert\xi} := N^{-1}\sum_{i=1}^{N}\bigl(x_{t\vert\xi}^{i}\bigr)^2\;,\text{ and}\quad
  \check{x}_{t\vert\xi} :=  N^{-1}\sum_{i=1}^{N}\bigl(x_{t\vert\xi}^{i}\bigr)^3\;.
\end{align*}
Thus, the quantities defining the matrices and
right-hand sides involved in the estimation step (cf.\
\eqref{eq:fun:form:param:sys} and \eqref{eq:fun:form:param:sys2},
respectively) are approximated via
\begin{subequations}
\begin{align}
  a_1(\xi)^T&\approx \frac{h}{2}\bigl(Q_n(\bar{x}_{\cdot\vert\xi}),
  Q_n(\check{x}_{\cdot\vert\xi})\bigr)^T\;,\quad
  b_1(\xi)\approx \bar{x}_{t\vert\xi} - \xi \label{eq:ex:drift:est}\\
  a_2(\xi)^T&\approx \frac{h}{2}\bigl(2n,
  Q_n(\tilde{x}_{\cdot\vert\xi})\bigr)^T\;,\quad
  b_2(\xi)\approx N^{-1}\sum_{i=1}^{N}
  \Bigl(
  x_{t\vert\xi}^{i} - \xi -
  Q_n\bigl(f_{\hat{\vartheta}}\circ x_{\cdot\vert\xi}^{i}\bigr)
  \Bigr)^2\;,\label{eq:ex:diffusion:est}
\end{align}
\label{eq:ex:est}%
\end{subequations}
with $f_{\hat{\vartheta}}\equiv f(\cdot;\hat{\vartheta})$ and
$Q_n$ denoting the quadrature operator of the trapezoidal rule on $[0,t]$ with
$n$ equally spaced ($h=t/n$) subintervals (cf.\ \eqref{eq:intapprox:trapez})
\begin{align}
  Q_n(u):=\frac{h}{2}(u_0 + u_{t} + 2\sum_{j=1}^{n-1}u_{jh})\;.\label{eq:trapez:op}
\end{align}

It should be emphasized that equations such as $x=A^{+}b$ (e.g.\ as in
lines $10$ and $22$ in Algorithm \ref{alg:est}) are merely meant as
a formal notation for $x$ solving the least squares problem
$x=\argmin_{s\in\mathcal{S}}\norm{s}_2^2$, $\mathcal{S} =
\{z\colon\norm{Az-b}_2^2\rightarrow\textrm{min}\}$, cf.\ equations
\eqref{eq:drift:est} and \eqref{eq:diffusion:est}, rather than
indicating that we first compute $A^{+}$ and then multiply it with
$b$ to obtain $x$, a step which clearly is computationally
inefficient. In practice the method to compute the solution
typically depends on the rank of $A$. For the numerical examples
that we will present in Section \ref{sec:num} we ensured a full
rank situation by considering a sufficiently large number of
different initial conditions. The solution $x$ can be computed 
with one of the several methodologies for solving least squares problems,
e.g.\ the Cholesky factorization of the normal equations\footnote{This method 
    for solving the least squares problem is simple to implement but might not be optimal 
    and is usually only recommended for problems with large residuals. However, it is 
    noteworthy that for the numerical examples presented in this study, the method 
    performed well and no stability issues occurred (as can be demonstrated by monitoring 
    the condition number; not shown here), and similar results were obtained with a 
    QR factorization with column pivoting (not shown here).} $A^TAx=A^Tb$; 
for details on such methodologies we refer to standard textbooks on
numerical linear algebra, e.g.\ \cite[ch.\ $5$]{Golub1996}.

\subsection{Properties of the Estimator}
\label{sec:estimator:properties}

The proposed estimation procedure relies on two key ingredients. The
first is that the methodology is based on the identities in
\eqref{eq:est:start}. The second is that by considering a finite
sequence of initial conditions we can cope with multiple parameters in drift
and diffusion coefficient, respectively. In the sequel we demonstrate the
influence of both components on the proposed estimation scheme
with the help of some elementary, nonetheless illustrative, examples
when no multiscale effects are present. A detailed and rigorous
analysis of the proposed methodology for both multiscale and
non-multiscale situations will be presented elsewhere.

To illustrate the influence on the identities in
\eqref{eq:est:start} and to address some asymptotic properties we
consider a simple Langevin equation with additive noise
\begin{align*}
dx_t = \vartheta f(x_t)\,dt + \sqrt{\theta}\,dW_t\;.
\end{align*}
The drift estimator proposed in this study -- this parameterization is
already a straightforward generalization of the one introduced above --
relies on the relation
\begin{align}
  \E\bigl(\int_0^tf(x_s)\,ds\bigr)\vartheta = \E(x_t-x_0)\;,
  \label{eq:1par:drift:martingale}
\end{align}
with $\vartheta$ being the true value. For a fixed final time
$t<\infty$, the estimator for continuous-time observations --
meaning that we approximate only the expectation by an average
  but do not approximate the integrals -- based on a single (fixed) initial condition $x_0=\xi$ is given by
\begin{align}
  \hat{\vartheta} &= \frac{\sum_{i=1}^N(x_t^{i}-x_0)}{\sum_{i=1}^N
    \int_0^tf\bigl(x_s^{i}\bigr)\,ds}
  = \vartheta + \frac{\sqrt{\theta}\sum_{i=1}^N\int_0^t\,dW_s^i}{\sum_{i=1}^N\int_0^tf(x_s^i)\,ds}
  =\vartheta + \frac{\mathcal{N}(0,\theta t/N)}{\frac{1}{N}\sum_{i=1}^N\int_0^tf(x_s^i)\,ds}\;.
  \label{eq:1par:drift:err:cont}
\end{align}
Notice that we dropped the dependency 
on the initial condition,
because only a single initial condition is considered here. Since we
approximate only the expectations by finite averages, it is not
surprising that the
estimator for continuous-time observations \eqref{eq:1par:drift:err:cont}
converges to the true value in agreement
with the law of large numbers. The property that the variance of the
error vanishes for $N\rightarrow\infty$ reflects the fact that the
estimator relies on an identity, i.e.\ on a direct (deterministic)
computation \eqref{eq:1par:drift:martingale} rather than on
asymptotic time limits (i.e.\ on ergodicity).

Recall that we have $(n+1)$ observations at times $0=t_0<t_1<\cdots < t_n=t$ with $t>0$
fixed, in the case of discrete-time observations. The integral
approximation introduces an additional error that can be identified
via
\begin{align}
  \int_0^t u_s\;ds = Q_n(u) + c_n\;, \label{eq:intapprox:error:worstcase}
\end{align}
with an appropriate constant $c_n\in\R$ that depends not only on $n$
but also on $u$ and $t$, and which vanishes in the limit as
$n\rightarrow\infty$. Here $Q_n$ denotes the quadrature operator of
the trapezoidal rule as defined in \eqref{eq:trapez:op}. The actual
error of the trapezoidal rule depends on the regularity of the
integrand. Here we rely only on the assumption $\lim_{n\rightarrow\infty}c_n=0$ as
a general scenario and do not discuss the rate of convergence. Then,
similarly to the continuous-time case, the estimator (when
neglecting sampling errors) can be written as
\begin{align}
  \hat{\vartheta} &= \frac{\sum_{i=1}^N(x_t^i-x_0)}{\sum_{i=1}^N Q_n(f\circ x_{\cdot}^i)}
  = \frac{\sum_{i=1}^N\bigl(\vartheta Q_n(f\circ x_\cdot^i) + \vartheta c_n^i +
    \sqrt{\theta}\int_0^t\,dW_s^i\bigr)}{\sum_{i=1}^N Q_n(f\circ x_\cdot^i)}
  \nonumber \\
  &= \vartheta \biggl(1-\frac{\bar{c}_N(n)}{\bar{c}_N(n)-\frac{1}{N}\sum_{i=1}^N\int_0^tf(x_s^i)\,ds}\biggr) +
  \frac{\mathcal{N}(0,\theta t/N)}{\frac{1}{N}\sum_{i=1}^N\int_0^tf(x_s^i)\,ds - \bar{c}_N(n)}\;,
  \label{eq:1par:drift:err:discrete}
\end{align}
with $\bar{c}_N(n) = N^{-1}\sum_{i=1}^Nc_n^i$ being the average
error constant. Hence, the integral approximation introduces an
additional bias -- second term in the bracket -- that can be
controlled by $n$. Notice that $\abs{\bar{c}_N(n)}\le \max_{1\le
i\le N}\abs{c_n^i}$ so that the additional bias vanishes as
$n\rightarrow\infty$ for every $N$. Consequently, from equation
\eqref{eq:1par:drift:err:discrete} one infers that for a fixed final
time $t >0$, one should choose both $n,N\gg 1$ to obtain an accurate
estimate.

To estimate the diffusion coefficient the proposed scheme relies on
the relation
\begin{align*}
  \theta = \frac{1}{t}\E\Bigl(\bigl(x_t-x_0-\vartheta\int_0^tf(x_s)\,ds\bigr)^2\Bigr)\;,
\end{align*}
that is valid for all times $t>0$ and where we replace $\vartheta$
by its estimator $\hat{\vartheta}$ for concreteness. Consequently,
the estimator for continuous-time observation using a single (fixed)
initial condition and a fixed final time $t>0$ reads
\begin{align*}
  \hat{\theta} &=\frac{1}{tN}\sum_{i=1}^N\Bigl(x_t^i-x_0-\hat{\vartheta}\int_0^tf(x_s^i)\,ds\Bigr)^2
  = \frac{1}{tN}\sum_{i=1}^N\Bigl((\vartheta-\hat\vartheta)\int_0^tf(x_s^i)\,ds
  + \sqrt{\theta}\int_0^t\,dW_s^i\Bigr)^2\\
  &= \theta\frac{1}{N}\chi_N^2
  + \frac{(\vartheta-\hat\vartheta)^2}{tN}\sum_{i=1}^N\Bigl(\int_0^tf(x_s^i)\,ds\Bigr)^2
  + \frac{2(\vartheta-\hat\vartheta)}{tN}\sum_{i=1}^NW_t^i\int_0^tf(x_s^i)\,ds\;.
\end{align*}
Since the estimator for the diffusion coefficient depends on the estimated drift parameter, an
additional error is introduced (last two terms), as expected. To illustrate the asymptotic
properties of the estimator we, however, assume that the error $(\vartheta-\hat\vartheta)$
is negligible. Consequently, we find that
\begin{align*}
  \hat{\theta} \approx \theta\frac{1}{N}\chi_N^2\;,
\end{align*}
where $\chi_N^2$ denotes the Chi-squared distribution with $N$ degrees of freedom.
Recall that $\frac{1}{N}\chi_N^2\approx \mathcal{N}(1,2/N)$ for $N$ sufficiently large, as a
consequence of the central limit theorem.

For the modification of discrete-time observations the integrals are again approximated via the
trapezoidal rule. Based on the same $(n+1)$ observations (neglecting sampling errors)
at $0=t_0<t_1<\cdots < t_n=t$ and $t>0$ fixed we find
\begin{align*}
  \hat{\theta} &=\frac{1}{tN}\sum_{i=1}^N\Bigl(x_t^i-x_0-\hat{\vartheta}Q_n(f\circ x_\cdot^i)\Bigr)^2
  = \frac{1}{tN}\sum_{i=1}^N\Bigl((\vartheta-\hat\vartheta)\int_0^tf(x_s^i)\,ds
  + \sqrt{\theta}\int_0^t\,dW_s^i + \hat\vartheta c_n^i\Bigr)^2\;,
\end{align*}
with $c_n^i$ being the error representation of the trapezoidal rule, cf.\
\eqref{eq:intapprox:error:worstcase}. In contrast to the continuous-time situation,
the term $\hat\vartheta c_n^i$ reflects the additional error due to the integral approximation
that can be controlled by $n$. Since this additional term is the only difference, expanding the
square yields a similar result as in the continuous-time situation.
If we assume again that the error from the drift parameter $(\vartheta-\hat\vartheta)$
is negligible, then we find
\begin{align}
  \hat{\theta} &\approx \theta\frac{1}{N}\chi_N^2 + \frac{{\hat\vartheta}^2}{t}\tilde{c}_N(n)
  + \frac{2\hat\vartheta\sqrt{\theta}}{tN}\sum_{i=1}^Nc_n^iW_t^i\;,
\label{eq:1par:diff:err:discrete}
\end{align}
with $\tilde{c}_N(n) = N^{-1}\sum_{i=1}^N{(c_n^i)}^2$. Since $0\le
\tilde{c}_N(n)$ for every $N$, the second term in equation
\eqref{eq:1par:diff:err:discrete} can only be controlled by $n$.
Consequently, for a fixed final time $t$,
choosing $n,N\gg 1$ is necessary to obtain an accurate estimate
of the true parameter. It is noteworthy that the same steps may be
carried out for an arbitrary diffusion function, but obviously we
cannot directly infer the distributions of terms involving the
stochastic integral.

To deal with multiple parameters in drift and/or diffusion
coefficients the proposed estimation scheme relies on considering a
finite sequence of initial conditions and defining the estimator via
the best approximation. To illustrate the effect of this second
component of the estimation procedure, consider the SDE
\begin{align*}
  dx_t = \vartheta f(x_t)\,dt + \sqrt{g(x_t)}\,dW_t.
\end{align*}
Provided all quantities are well defined, the estimator of
$\vartheta$, via a best approximation using a sequence of initial
conditions $(\xi_{i})_{1\le i\le m}$, is given by
\begin{align*}
\hat{\vartheta} = \frac{\sum_{i=1}^ma_1(\xi_{i})b_1(\xi_{i})}{\sum_{i=1}^ma_1(\xi_{i})^2}\;,
\end{align*}
where $a_1(\cdot)$ and $b_1(\cdot)$ are as in
\eqref{eq:fun:form:param:single} associated with the drift function
$f$. On the other hand, the quantity $b_1(\xi_{i})/a_1(\xi_{i})$
yields a local estimator for each initial condition because the
considered problem has only one parameter to be determined. Thus the
best approximation corresponds here to the weighted arithmetic mean
of these local estimators with weights $a(\xi_{i})^2$. Consequently,
increasing $m$ includes an additional stabilization effect into the
estimation scheme. From this point of view, the best approximation
resolves naturally the problem of combining local estimates to a
global estimator that arises in different estimation procedures as
well, where piecewise local strategies are utilized to improve the
estimates; see for example \cite{Calderon2007}. In this study a
heuristic estimation strategy for non-constant diffusions in the MLE
framework is proposed by extracting local information from a large
time-series to estimate coefficients locally and combine these local
estimators to global estimators. Since this strategy is based on the
MLE, the results in the above study also highlight the subsampling
issue when applied to multiscale diffusions.

It should be noted that although the proposed methodology computes
moments of the solution of an SDE, there is no direct link to the
generalized method of moments (GMM) \cite{Hansen1982}. 
The GMM for parametric estimation for SDEs is more closely related to
the MLE instead, as both estimation schemes rely on ergodicity of
the corresponding process by using long time-series. In fact, 
both estimators even coincide for many cases 
\cite[ch.\ $14.4$]{Hamilton1994}.
Since the MLE is biased for multiscale diffusions, it can be expected that 
the subsampling issue also arises for the GMM
when applied to data obtained form a system with multiscale structure.

\subsection{Estimators for Multiscale Diffusions}
\label{sec:estimator:multiscale}
In the context of diffusion processes with two widely separated time scales we consider the following set-up: A fast/slow system of SDEs
\begin{subequations}
  \begin{align}
    dx_t &= \Bigl(\frac{1}{\varepsilon}f_1(x_t,y_t) + f_0(x_t,y_t)\Bigr)\,dt + \alpha(x_t,y_t) \,dU_t\;,\label{eq:fastslow:slow}\\
    dy_t &= \Bigl(\frac{1}{\varepsilon^2}g_2(x_t,y_t) + \frac{1}{\varepsilon}g_1(x_t,y_t) + g_0(x_t,y_t)\Bigr)\,dt
    + \frac{1}{\varepsilon}\beta(x_t,y_t) \,dV_t\;,\label{eq:fastslow:fast}
  \end{align}
  \label{eq:fastslow}%
\end{subequations}%
equipped with appropriate initial conditions. In \eqref{eq:fastslow}
$U,V$ denote Brownian motions of appropriate dimensions and
$0<\varepsilon\ll 1$ denotes a small parameter. For the dimension of
the fast/slow system we assume that $y\colon T\mapsto\R^d$ and (for
simplicity) $x\colon T\mapsto\R$, where $T=[0,t]$ denotes a finite
time interval of interest. Furthermore, we assume that drift and
diffusion functions in \eqref{eq:fastslow:slow} and
\eqref{eq:fastslow:fast} respectively, are such that there exists a
well-defined (i.e. the SDE provides a unique weak solution on any
finite time interval and for any initial condition) coarse-grained
SDE
\begin{align}
  dX_t = f(X_t)\, dt + \sqrt{g(X_t)}\,dW_t\;,\label{eq:fastslow:ampl}
\end{align}
in the limit of $\varepsilon\rightarrow 0$; see e.g.\ \cite[ch.~11 and 18]{Pavliotis2008book} and
references therein for technical details. That is, the slow process $x$ is
approximated by the solution of \eqref{eq:fastslow:ampl} for $\varepsilon\ll 1$.
Even in cases where both
the effective drift $f$ and the effective diffusion function $g$ are known, the actual
computation of these expressions might be difficult or even impossible,
as appropriate
Poisson equations have to be solved and integrals with respect to the invariant measure
of the fast process to be computed.
Hence, our goal
is to estimate both effective drift coefficient $f$ and effective diffusion coefficient $g$ in
\eqref{eq:fastslow:ampl} from available data (observations) of the fast/slow system \eqref{eq:fastslow}
(more precisely, only of its slow component). To this end we assume the same parameterizations of the effective drift and
diffusion functions as introduced in Section \ref{sec:estimator:derivation}:
\begin{align*}
f(x)\equiv f(x;\vartheta) := \sum_{j\in J_f} \vartheta_j x^j\quad\text{and}\quad
g(x)\equiv g(x;\theta) := \sum_{j\in J_g} \theta_j x^j\;,
\end{align*}
where we recall that $J_f,J_g \subset \N_0$ denote index sets with
$p=\abs{J_f}$ and $q=\abs{J_g}$. Our goal then is: Given
observations only of the slow component \eqref{eq:fastslow:slow}, is it
possible to estimate the parameters
$\vartheta\equiv{(\vartheta_j)}_{j\in J_f}\in\R^{p}$ and
$\theta\equiv {(\theta_j)}_{j\in J_g}\in\R^{q}$ characterizing the associated
coarse-grained equation \eqref{eq:fastslow:ampl}? Under the
assumption that the (ergodic) fast process \eqref{eq:fastslow:fast}
is stationary, the algorithm described in Section
\ref{sec:estimator:derivation} applies straightforwardly also for
this problem, given the final time $t$ of observation length is
appropriately chosen; cf.\ Section \ref{sec:num:multi}. That is,
using multiple initial conditions for the slow process
\eqref{eq:fastslow:slow} to deal with multiple parameters in drift
and diffusion, while the fast process is sampled from its invariant
measure that is assumed to be known, either analytically or
numerically.

The main motivation to consider the same algorithm also in the
presence of multiscale effects originates from the fact that both
the slow component of the full fast/slow system and the effective
dynamics have probability laws that are approximately the same, 
provided $\varepsilon\ll 1$.
Consequently, also expectations with respect to these laws are
approximately equal. Since the proposed methodology is based
on expectations, cf.\ equation \eqref{eq:est:start},  
we believe that this approach yields
asymptotically unbiased estimators. A rigorous analysis of the
proposed methodology to verify this heuristic argument is currently
work in progress and will be reported elsewhere.

We conclude this section with a remark on a recently proposed
estimator for constant diffusion coefficients \cite{Frederix2011a}
that can also be derived using the approach we introduce here. 
To this end, assume that the coarse-grained
equation takes the form:
\begin{align*}
  dx_t =  f(x_t)\,dt + \sqrt{\theta}\, dW_t\;,
\end{align*}
where we assume that the effective drift $f$ is known and we wish to estimate the diffusion coefficient
$\theta$ from available data of a fast/slow system. When considering a single (fixed) initial condition
and following the approach introduced in Section \ref{sec:estimator:derivation}, the resulting estimator
reads
\begin{align*}
  \hat{\theta} =
  \frac{1}{tN}\sum_{i=1}^{N}\Bigl(\bigl(x_t^i - x_0-\sum_{j=0}^{n-1}\int_{jh}^{(j+1)h}f(x_s^i)\, ds\bigr)^2\Bigr)\;,
  \end{align*}
where the integrals are being approximated by a quadrature rule and
$t$ is chosen appropriately (in fact, we have approximated the
Lebesgue integral via the trapezoidal rule). If one, however, uses a
rectangular-method with the left corner node instead, this estimator
coincides with the estimator proposed in \cite{Frederix2011a} for
estimating the effective diffusion coefficient based on observations
of the slow component of a fast/slow system.
We emphasize, that a crucial assumption on
the estimator in the aforementioned work is that the effective drift
is known a priori. This assumption is very restrictive and makes the
estimator unfeasible for most practical applications. A further
limitation of the estimator in \cite{Frederix2011a} is that it
applies only in situations where the noise in the coarse-grained
equation is additive (constant diffusion coefficient). Conversely,
the methodology proposed here aims to estimate multiple parameters
in both the drift and the diffusion coefficients.


%% file: numerics.tex
\section{Numerical Experiments}
\label{sec:num} We now present numerical experiments of parameter
estimations for diffusion processes to illustrate the behavior of
the estimation scheme developed in Section \ref{sec:estimator}.

\subsection{Parameter Estimation for single-scale SDEs}
\label{sec:num:single} Let us first present numerical results for
parameter estimation when no multiscale effects are present. We
investigate two different SDEs. In Section \ref{sec:num:single:OU}
we consider the SDE corresponding to the Ornstein-Uhlenbeck process
and in Section \ref{sec:num:single:GL} we examine the stochastic Landau-Stuart
equation. The purpose of these examples is twofold. On the one
hand they are used to illustrate that the proposed methodology can
be employed successfully in practice to estimate unknown parameters
and on the other they help us understand the influence of the
parameters $n$, $h$, $N$, and $m$ on the algorithmic estimation
procedure. The time series were obtained by solving the
corresponding SDEs via the Euler-Maruyama scheme. In all numerical
examples reported in this section, we used a time step of
$h=10^{-3}$.

\subsubsection{Ornstein-Uhlenbeck Process}
\label{sec:num:single:OU}
Consider the following SDE
\begin{align}
dx_t = -A x_t \,dt + \sqrt{\sigma}\,dW_t\;,\quad x_0=\xi\;,\label{eq:num:OU}
\end{align}
with the unique solution being the Ornstein-Uhlenbeck process starting at $\xi$.
The estimation procedure is applied to data
from the SDE with true parameters $(A,\sigma) = (0.5,0.5)$ using a variety of different values
for the parameters of the algorithm.
\begin{figure}[]
  \centering
  \subfigure[relative error of $\hat{A}$]{
    \includegraphics[width=0.465\textwidth]{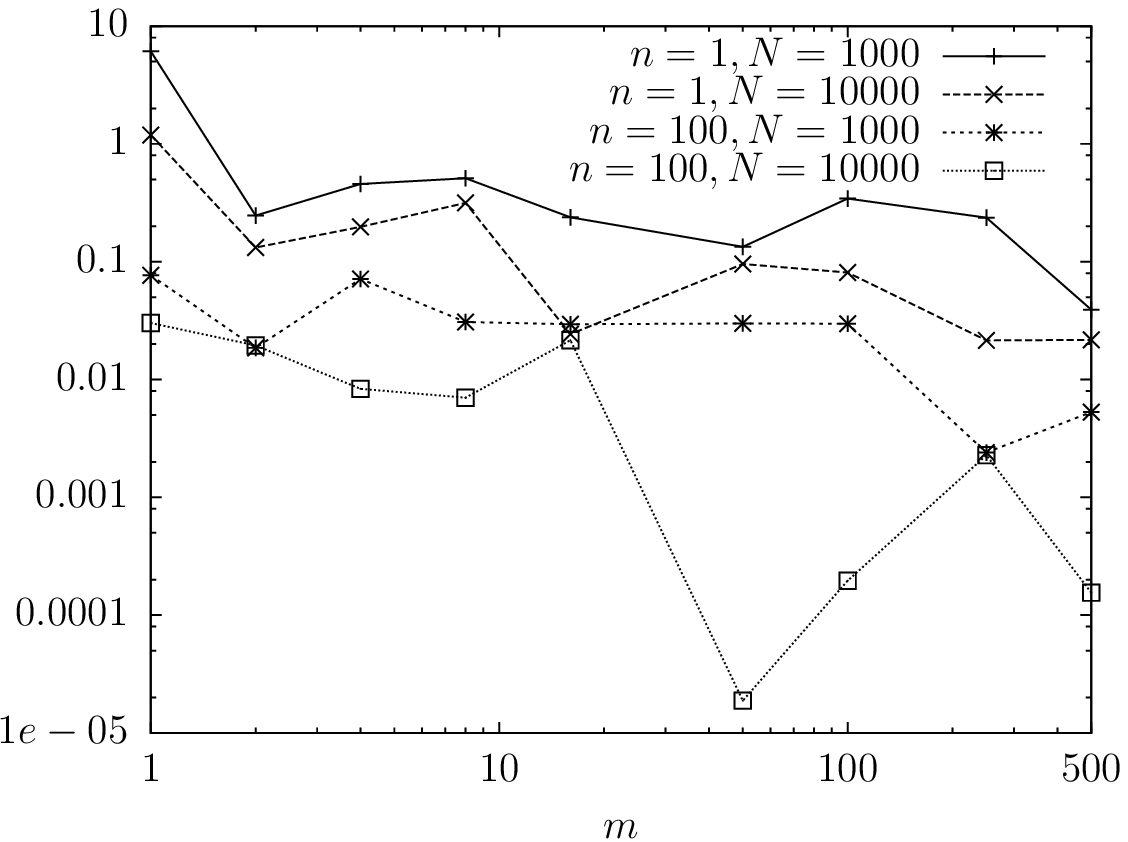}
    \label{fig:nosub:ou:p1}
  }
  \quad
  \subfigure[relative error of $\hat{\sigma}$]{
    \includegraphics[width=0.465\textwidth]{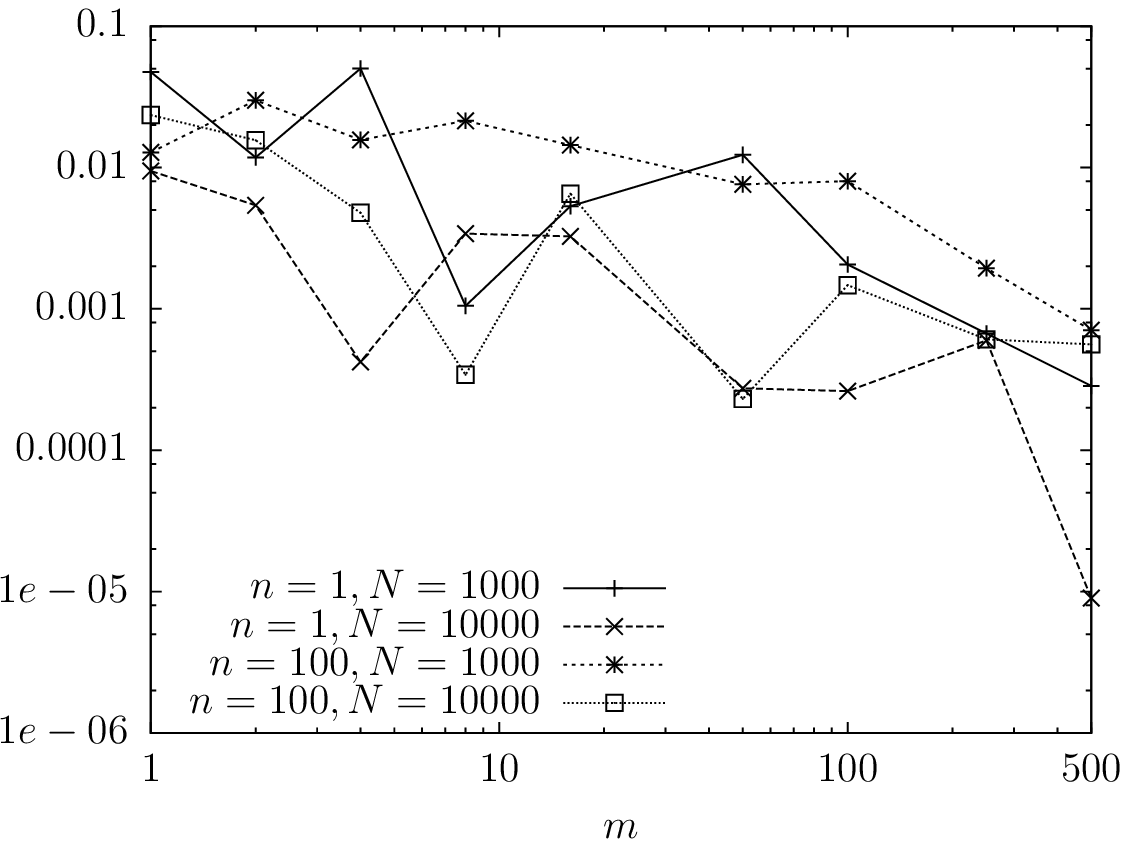}
    \label{fig:nosub:ou:p2}
  }
  \caption[]{Relative error of the estimated parameters in \eqref{eq:num:OU} as functions of the
    number of initial conditions $m$ in a log-log scale. The final time of
    the considered time series is $t=nh$ with $h=0.001$ and the true parameters are
    $(A,\sigma) = (0.5,0.5)$. Furthermore, $N$ denotes the number of
    independent Brownian paths.}
  \label{figure:nosub:ou}
\end{figure}
Figure \ref{figure:nosub:ou} depicts the relative errors of the
estimated values as functions of the number of initial conditions
$m$ for different combinations of $n$
(recall that the final time is $t=nh$) and $N$.
For the estimated drift parameter
$\hat{A}$ (Figure \ref{fig:nosub:ou:p1}) there is a discrepancy
among different combinations of $n$ and $N$ for small values of $m$
and the relative errors vary from approximately $0.08$ to $6$.
Increasing $m$ decreases the relative errors, with fluctuations due
to the discretization, as expected; cf.\ Section
\ref{sec:estimator:properties}. Moreover, it is apparent that the
larger $n$ and $N$, the smaller the relative error. In contrast to
the drift parameter, the relative errors of the diffusion parameter
$\hat{\sigma}$ (Figure \ref{fig:nosub:ou:p2}) are already small
(relative error smaller than $0.05$) even for small values of $m$. This is due to the
constant diffusion coefficient of the SDE. Increasing $m$ generally
decreases the relative error further, again, with fluctuations due
to the discretization. We note that although an error propagates from
the drift estimation to the diffusion estimation (cf.\ Section
\ref{sec:estimator}), it appears to be negligible in this example.
Based on theses results, it seems plausible to tune the
algorithm-defining parameters in a way such that the estimators
provide a given relative accuracy while at the same time the
computational cost is minimized. Obviously the question of optimized
algorithm-defining parameters is of high importance in practical
applications. However, this is not a straightforward issue to
address as not only the discretization related parameters $n$ and
$N$ influence the accuracy of the algorithm in practice, but also
the number of initial conditions $m$ as well as their locations; we
will investigate these and related issues in a future study, see also the
discussion in Section \ref{sec:outlook}.

\subsubsection{Landau-Stuart Equation}
\label{sec:num:single:GL} Consider the stochastic Landau-Stuart
equation \cite[ch.\ $2.2$]{Kuramoto2003}, where both additive and
multiplicative noise are present
\begin{align}
dx_t = \bigl(Ax_t - Bx_t^3\bigr)\,dt + \sqrt{\sigma_a + \sigma_bx_t^2}\,dW_t\;,\quad
x_0=\xi\;.\label{eq:num:GL}
\end{align}
This SDE can be obtained from a wide class of spatially extended
systems, e.g. the noisy Kuramoto-Sivashinsky
equation~\cite{Pradas2011,Pradas2012} by assuming near-critical
conditions, i.e.\ being sufficiently close to the primary
bifurcation, and employing the homogenization theory developed
in~\cite{Blomker2007}.
In this case we need to estimate a total number of four parameters,
two in the drift and two in the diffusion coefficient.
\begin{figure}[]
  \centering
  \subfigure[relative error of $\hat{A}$]{
    \includegraphics[width=0.465\textwidth]{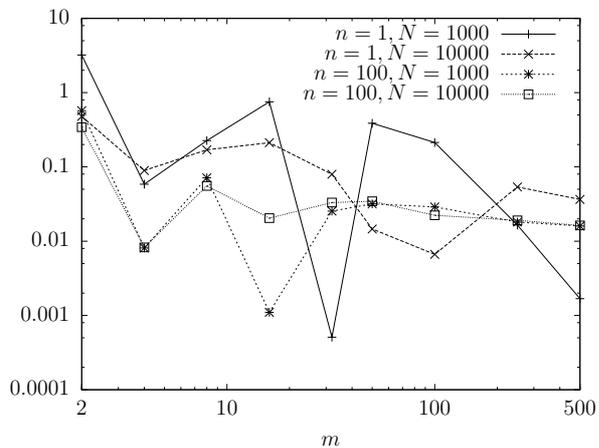}
    \label{fig:nosub:ls:p1}
  }
  \quad
  \subfigure[relative error of $\hat{B}$]{
    \includegraphics[width=0.465\textwidth]{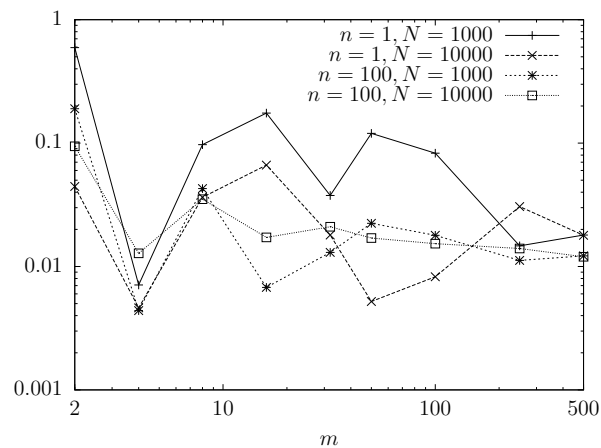}
    \label{fig:nosub:ls:p2}
  }
  \quad
  \subfigure[relative error of $\hat{\sigma}_a$]{
    \includegraphics[width=0.465\textwidth]{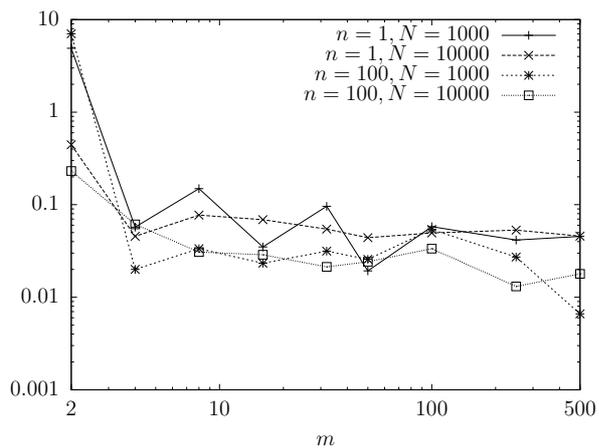}
    \label{fig:nosub:ls:p3}
  }
  \quad
  \subfigure[relative error of $\hat{\sigma}_b$]{
    \includegraphics[width=0.465\textwidth]{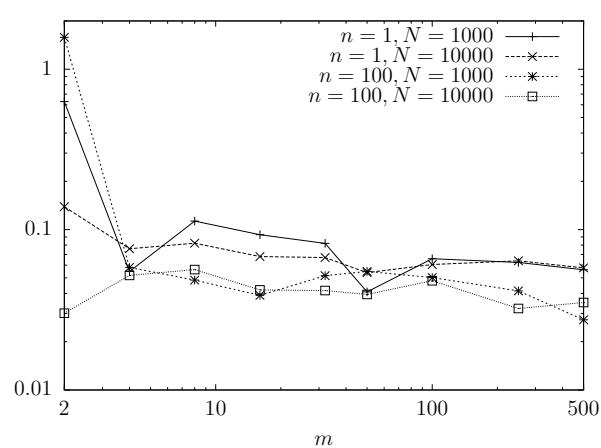}
    \label{fig:nosub:ls:p4}
  }
  \caption[]{Relative error of the estimated parameters in \eqref{eq:num:GL} as functions of the
    number of
    initial conditions $m$ using a log-log scale. The final time of
    the considered time series is $t=nh$ with $h=0.001$ and the true parameters are
    $(A,B,\sigma_a,\sigma_b) = (3,2,1.5,1.3)$. $N$ denotes the number of
    independent Brownian paths.}
  \label{figure:nosub:ls}
\end{figure}
We performed various numerical experiments with different choices of
parameters.
Figure
\ref{figure:nosub:ls} illustrates the relative error of the
estimated parameters as functions of the number of initial
conditions for different combinations of $n$ and $N$ when the true
parameters are $(A,B,\sigma_a,\sigma_b) = (3,2,1.5,1.3)$. We find
qualitatively the same behavior as in the previous section for the
Ornstein-Uhlenbeck process: all three parameters
$n$, $N$, and $m$ affect the accuracy of the estimators. Although the Figures
\ref{fig:nosub:ls:p1}--\subref{fig:nosub:ls:p4} show a decreasing
trend of the relative errors when increasing $m$, fluctuations are
still present. These fluctuations are of different magnitude for
different parameters and are reduced by increasing both $n$ and $N$.
One also observes that increasing $m$ improves the accuracy of the
estimators but only up to a certain level, i.e.\ the corresponding
curves approach nearly constant values for large $m$, that are due
to the approximation of the original system matrix, for instance in
\eqref{eq:fun:form:param:sys}, by a matrix based on
observations and the aforementioned discretizations. Although it is
apparent that estimating all parameters in the Landau-Stuart model
is more delicate than for the Ornstein-Uhlenbeck process in Section
\ref{sec:num:single:OU}, it nonetheless seems possible to determine
optimal algorithm-defining parameters such that the computational
cost is minimized given a certain error tolerance for the estimators
for this model also.

\subsection{Parameter Estimation for Fast-Slow Systems}
\label{sec:num:multi} Of particular interest is the behavior of the
estimator when applied to systems with two different time scales. We
examine the properties of the estimation scheme for stochastic
multiscale diffusions (Sections
\ref{sec:num:multi:fastOU1}--\ref{sec:num:multi:langevin2d}),
the problem of estimating the eddy diffusivity in a two-dimensional
  cellular flow (the Taylor-Green flow, Section \ref{sec:num:multi:taylorgreen}),
truncated systems of time rescaled
stochastic partial differential equations (Section
\ref{sec:num:multi:burgers}), and deterministic systems that exhibit
temporal chaos that can be approximated by an appropriate SDE
(Section \ref{sec:num:multi:chaos1}).
To measure the accuracy of the estimation procedure in these examples,
we rigorously derive the coarse-grained equations from the associated fast/slow
systems using homogenization theory so that
the theoretical coefficients are known. Based only on observations of the
slow component of the fast/slow system, the goal is to infer the coefficients
in the coarse-grained equation using the proposed estimation procedure. The
estimated values are then compared with the theoretical ones.
It is noteworthy that no
assumptions on the knowledge of the fast component, nor on the structure
of the fast/slow system are made (also $\varepsilon$ is unknown).
As the precise dependency of the estimation procedure on the control
parameters $n,m,N,h$ is still an open question, the main purpose of
this section is to illustrate the general applicability of the
proposed estimation procedure for multiscale diffusions. If not
stated otherwise, the generated time series were obtained by solving
the corresponding multiscale SDEs via the Euler-Maruyama scheme
using a time step $h=10^{-3}$. Furthermore, the expectation is
approximated by an average using $N=5000$ independent Brownian paths and
$m=150$ different (equally-spaced) initial
conditions are used. We emphasize once again that this particular
choice of algorithm-defining parameters might be far from optimal in
the sense of computational complexity.
But since our main goal is to demonstrate the applicability of the
proposed scheme to multiscale diffusions, these algorithm-defining
parameters will yield reliable estimators.

\subsubsection{Fast Ornstein-Uhlenbeck Noise}
\label{sec:num:multi:fastOU1} When the fast process is an
Ornstein-Uhlenbeck process it is rather straightforward to determine
the precise form of the effective equation associated to the
fast/slow system, because this task reduces to computing Gaussian
integrals. Consider for example
\begin{subequations}
\begin{align}
  dx_t &= \Bigl(\frac{1}{\varepsilon}\sigma(x_t)y_t + h(x_t,y_t) - \sigma'(x_t)\sigma(x_t)\Bigr)\,dt\;,
  \label{eq:sub:OU:fastslow:slow}\\
  dy_t &= -\frac{1}{\varepsilon^2}y_t\,dt + \frac{\sqrt{2}}{\varepsilon}\,dV_t\;,\label{eq:sub:OU:fastslow:fast}
\end{align}
\label{eq:sub:OU:fastslow}%
\end{subequations}
with $V_t$ being a standard Brownian motion, then the effective dynamics is given by
\begin{align}
  dX_t = \bar{h}(X_t)\,dt + \sqrt{2\sigma(X_t)^2}\,dW_t\;,\label{eq:sub:OU:effective}
\end{align}
where $\bar{h}(x)$ denotes the average of $h(x,\cdot)$ with respect
to the invariant measure of the fast process (Ornstein-Uhlenbeck
process). We note that we have subtracted the Stratonovich
correction from the drift in \eqref{eq:sub:OU:fastslow:slow}, so
that the noise in \eqref{eq:sub:OU:effective} can be interpreted in
the It\^o sense.\footnote{The noise
  entering \eqref{eq:sub:OU:fastslow:slow}, i.e.\ the process
  \eqref{eq:sub:OU:fastslow:fast}, is a smoothed approximation to white noise, so that the noise
  in the limiting equation has to be interpreted in the Stratonovich sense according to the
  Wong-Zakai theorem. Correcting the drift is not essential for the applicability of our
  methodology, it was done so that the limiting equation \eqref{eq:sub:OU:effective} is
  somewhat simpler.}
  In the sequel we consider two different choices of the pair $h(\cdot ),\sigma(\cdot )$.
  As a first example let
\begin{align}
  h(x,y) = h(x) = Ax\quad\text{and}\quad\sigma(x) = \sqrt{\sigma}\;,\label{eq:sub:OU:effective:OU}
\end{align}
then the amplitude equation is
precisely the Ornstein-Uhlenbeck process (cf.\ Section \ref{sec:num:single:OU}) but here in the context of multiscale
diffusions where classical estimators fail. To illustrate this failure and motivate the necessity of an
appropriate sampling rate, the quadratic variation of the path (QVP) estimator for the effective diffusion
constant and the maximum likelihood estimator (MLE) for the effective drift parameter are applied
to a time series on $[0,5000]$ with initial condition $x_0=0.5$ generated by the
associated fast/slow system \eqref{eq:sub:OU:fastslow} with true parameters
$(A,\sigma) = (-0.5,0.5)$ and $\varepsilon=0.1$.
\begin{figure}[]
  \centering
  \subfigure[Classical estimators]{
    \includegraphics[width=0.465\textwidth]{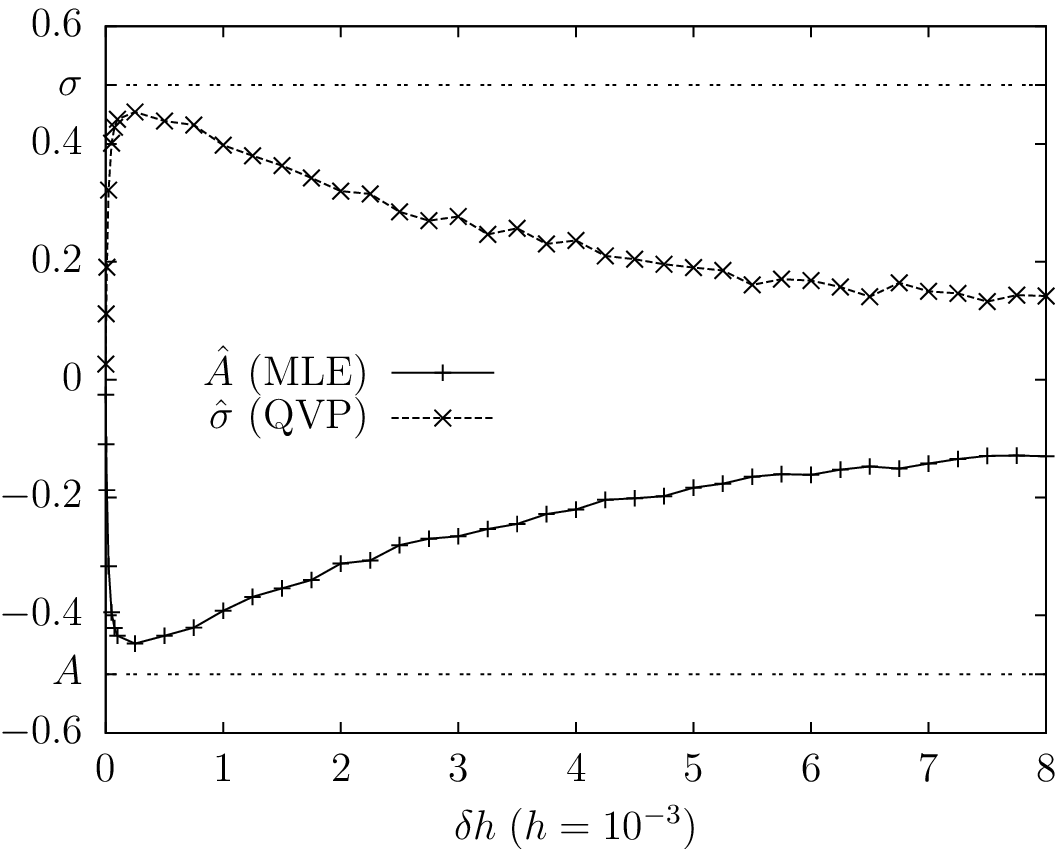}
    \label{fig:sub:ou:p1}
  }
  \quad
  \subfigure[Novel procedure]{
    \includegraphics[width=0.465\textwidth]{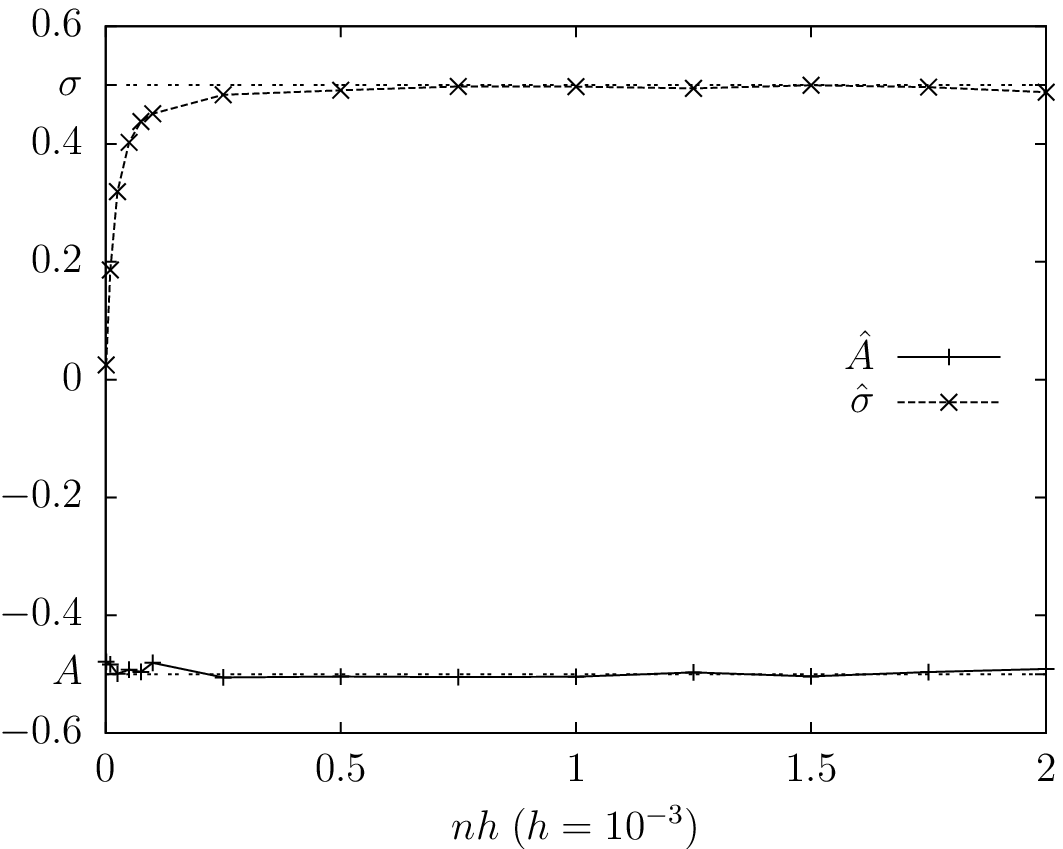}
    \label{fig:sub:ou:p2}
  }
  \caption[]{Performance of classical estimators and our procedure for both
    drift and diffusion coefficients in \eqref{eq:sub:OU:effective:OU}. In Figure
    \subref{fig:sub:ou:p1} the MLE and QVP as functions of the subsampling rate
    $\delta h$ ($\delta =1$ corresponds to no subsampling) are shown. Figure
    \subref{fig:sub:ou:p2} shows the result of our procedure as a function of the final time
    $t=nh$. In both cases the sampling rate of the considered time series is $h=0.001$
    and the true parameters are $(A,\sigma) = (-0.5,0.5)$. Dashed lines denote the true values.}
  \label{figure:sub:ou}
\end{figure}
Figure \ref{fig:sub:ou:p1} depicts the performance of both the
QVP and the MLE as functions of the
subsampling rate $\delta h$. The parameter $\delta$ indicates here
that only every $\delta$-th observation is used to estimate the
parameters in the drift and diffusion coefficients, hence $\delta h$ denotes the
time period between two consecutive observations.
Starting from the case without subsampling
($\delta=1$, that is $\delta h = 10^{-3}$) and increasing $\delta$,
both estimators approach the true value. However, after an optimal
subsampling rate, for which the estimator is as close to the true
value as possible (here approximately $\delta h=0.25$, that is
$\delta = 250$), both estimators deviate monotonically from the
target value. We note that at the optimal subsampling rate the
relative error is approximately $10\%$.
Figure \ref{fig:sub:ou:p1} seems to suggest that the optimal subsampling rate
is given by the local extremum of the estimator as function of the subsampling.
However, this behavior is not true in general, see for
instance the numerical examples in \cite{Pavliotis2007,Frederix2011a} that
reveal different behaviors. Recall, that the optimal
subsampling rate is in general unknown.

The situation is different when the effective parameters are
estimated via the method introduced here. Observations are
generated by the associated fast/slow system
\eqref{eq:sub:OU:fastslow} with true parameters $(A,\sigma) =
(-0.5,0.5)$ and $\varepsilon=0.1$. The performance of the estimator
as function of the final time $t=nh$ for both $\hat{A}$ and
$\hat{\sigma}$ is plotted in Figure \ref{fig:sub:ou:p2} and is
compared directly to the results obtained by the classical
estimators discussed before. For small values of $t=nh$ one observes
that the estimated value of the drift parameter $\hat{A}$ fluctuates
around the true value and stabilizes for larger times with minor
fluctuations around the true value. We note that the estimator
obtained by the proposed scheme significantly outperforms the MLE in
terms of accuracy.\footnote{At the optimal subsampling rate the relative error for the MLE is
  approximately 10\%, whereas the relative error for the drift parameter obtained by the novel scheme is
  less than 1\% for $t=nh\ge 0.25$. Admittedly, this comparison is not completely fair, because
  the novel estimation scheme in its current form relies on more data than the MLE, but we believe that
  these numbers nonetheless illustrate the potential of the methodology we propose.}
For the estimated diffusion coefficient
$\hat{\sigma}$ one finds that the novel scheme proposed here and the
QVP estimator yield similar results for small values of $t=nh$,
indicating that increasing $t=nh$ here has the same beneficial
effect as subsampling does for the QVP estimator. But unlike the
QVP, the estimator proposed here approaches the correct value
further and closely fluctuates around it when increasing $nh$,
without any deviations as when using QVP and subsampling. This is a
typical behavior for the estimator when applied to multiscale
diffusions as we will see in the forthcoming examples.
Consequently, one finds that, unlike classical estimators, once the
final time $t=nh$ is larger than a critical value, the estimator
fluctuates closely around the true value.\\

Consider as a second example $h(x,y)=h(x)=Ax-Bx^3$ and
$\sigma(x) = \sqrt{\sigma_a+\sigma_bx^2}$, so that the fast/slow
system \eqref{eq:sub:OU:fastslow} reads
\begin{subequations}
  \begin{align}
    dx_t &= \Bigl(\frac{1}{\varepsilon}y_t\sqrt{\sigma_a+\sigma_b{x_t}^2}
    + (A-\sigma_b)x_t-B{x_t}^3\Bigr)\,dt\;,
    \label{eq:sub:OU:QVP:fastslow:slow}\\
    dy_t &= -\frac{1}{\varepsilon^2}y_t\,dt + \frac{\sqrt{2}}{\varepsilon}\,dV_t\;,
    \label{eq:sub:OU:QVP:fastslow:fast}
  \end{align}
  \label{eq:sub:OU:QVP:fastslow}%
\end{subequations}
with the effective dynamics in \eqref{eq:sub:OU:effective} given by
the Landau Stuart equation (see Section \ref{sec:num:single:GL})
\begin{align}
  dX_t = (AX_t-B{X_t}^3)\,dt + \sqrt{2(\sigma_a+\sigma_b{X_t}^2)}\,dW_t\;.\label{eq:sub:OU:QVP:effective}
\end{align}
A natural extension of the QVP estimator to diffusion coefficients
that depend on multiple parameters is obtained by considering the
standard QVP relation for different time increments. Provided that
one considers a sufficient number of increments, i.e.\ a sufficient number
of estimating equations,
one can define the QVP estimator by solving the
arising system. To illustrate that this approach can deal with
multiple parameters in the diffusion coefficient when no multiscale
effects are present in the data, we first use the
QVP to estimate $\sigma_a$ and $\sigma_b$ in
\eqref{eq:sub:OU:QVP:effective} based on data that are also obtained
from \eqref{eq:sub:OU:QVP:effective}. This (single-scale) situation
corresponds to the classical case of parametric estimation and the
QVP can indeed be used to obtain accurate
estimators $\hat\sigma_a$ and $\hat\sigma_b$, as it is illustrated
in table \ref{tab:qvp:ls4}\subref{tab:qvp:ls4:single}.
\begin{table}
  \centering
  \subtable[Classical Setting]{
    \begin{tabular}{lcc}
      \toprule
      $N$ & $\hat\sigma_a$ & $\hat\sigma_b$ \\
      \midrule
      $1$ & $0.948301$ & $0.251906$ \\
      $10$ & $0.871561$ & $0.390639$ \\
      $100$ &   $0.817719$ & $0.480469$ \\
      $1000$ & $0.806024$ & $0.500243$ \\
      \bottomrule
    \end{tabular}
    \label{tab:qvp:ls4:single}
  }\qquad
  \subtable[Multiscale Setting]{
    \begin{tabular}{lcc}
      \toprule
      $N$ & $\hat\sigma_a$ & $\hat\sigma_b$ \\
      \midrule
      $1$ & $0.056463$ & $0.000000$ \\
      $10$ & $0.039087$ & $0.028682$ \\
      $100$ &   $0.040613$ & $0.026214$ \\
      $1000$ & $0.040363$ & $0.026727$ \\
      \bottomrule
    \end{tabular}
    \label{tab:qvp:ls4:multi}
  }
  \caption{QVP estimators when fitting the Landau-Stuart SDE
    \eqref{eq:sub:OU:QVP:effective} to observed data. For the results outlined in table
    \subref{tab:qvp:ls4:single} the observed data were obtained from \eqref{eq:sub:OU:QVP:effective},
    i.e.\ the classical setting, whereas for the results in table \subref{tab:qvp:ls4:multi} observations
    of the slow component \eqref{eq:sub:OU:QVP:fastslow:slow} of the fast/slow system
    \eqref{eq:sub:OU:QVP:fastslow} were used, i.e.\ the multiscale setting. In both cases the true
    parameters are $(\sigma_a,\sigma_b)=(0.81,0.49)$ and $\varepsilon=0.1$ was used in the multiscale
    setting.
    The parameter $N$ indicates the number of independent Brownian paths that have been used to compute
    an average.}
  \label{tab:qvp:ls4}
\end{table}
Therein the obtained estimators $\hat\sigma_a$ and $\hat\sigma_b$
are presented for different values of $N$, where $N$ indicates
the number of independent Brownian paths that have been used to compute
an average to improve the accuracy. The
Landau-Stuart equation \eqref{eq:sub:OU:QVP:effective} was solved
numerically on $T=[0,1000]$ starting at $X(0) = 0.5$ with true
parameters $(A,B,\sigma_a,\sigma_b)=(1,2,0.81,0.49)$. Apparently,
the QVP yields very accurate estimates in this
setting. However, things change when the QVP is adopted in
the presence of multiple time scales. In this
multiscale setting we wish to estimate the parameters $\sigma_a$ and
$\sigma_b$ in the effective dynamics \eqref{eq:sub:OU:QVP:effective}
from observations of the slow component
\eqref{eq:sub:OU:QVP:fastslow:slow} of the fast/slow system
\eqref{eq:sub:OU:QVP:fastslow}. The same true parameters and
configuration (i.e.\ same number and length of time increments) of
the QVP as in the classical example without
multiscale effects were used, since the QVP performed well therein.
Table \ref{tab:qvp:ls4}\subref{tab:qvp:ls4:multi} displays the
obtained estimators  $\hat\sigma_a$ and $\hat\sigma_b$ for different
values of $N$ when the observations are obtained from
\eqref{eq:sub:OU:QVP:fastslow} with scale separation
$\varepsilon=0.1$. Increasing $N$ yields QVP estimators with minor
fluctuations but both estimators are strongly biased, as expected.
Hence, an appropriate subsampling of the data would be required to
remove the bias, but, once again the optimal subsampling rate is not
known a priori and might, as in the case of the MLE, even be different for different
parameters.

Conversely, we will use
\begin{figure}[]
  \centering
  \includegraphics[width=0.465\textwidth]{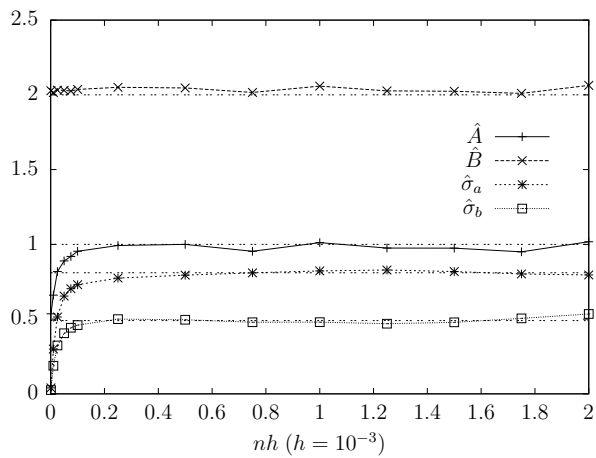}
  \caption[]{Performance of the novel estimators $\hat{A},\hat{B},\hat{\sigma}_a$, and
    $\hat{\sigma}_b$ for the Landau-Stuart equation \eqref{eq:sub:OU:QVP:effective} as functions
    of the final time $t=nh$ with $h=0.001$. The true effective parameters are
    $(A,B,\sigma_a,\sigma_b) = (1,2,0.81,0.49)$.}
  \label{figure:sub:gl}
\end{figure}
this example to illustrate that the parameters in multiscale
diffusions can be estimated accurately using the proposed scheme,
even though the amplitude equation provides a far more involved
structure than the one of the previous example in addition to the
multiscale structure of the problem. Figure \ref{figure:sub:gl}
depicts the performance of the estimation procedure based on
observations generated by the fast/slow system
\eqref{eq:sub:OU:QVP:fastslow} with true parameters $
(A,B,\sigma_a,\sigma_b)= (1,2,0.81,0.49)$ and $\varepsilon=0.1$ as a
function of the final time $t=nh$. The true values are indicated by
dashed lines. The behavior of the estimators is
qualitatively similar with that in the previous example. By increasing $t=nh$,
the estimators approach the true values, respectively and
fluctuate closely around them after a critical final time.
Consequently, all parameters can be estimated accurately.

\subsubsection{Brownian motion in a two-scale potential: A quadratic potential in one dimension}
\label{sec:num:multi:langevin}
\begin{figure}[]
  \centering
  \includegraphics[width=0.465\textwidth]{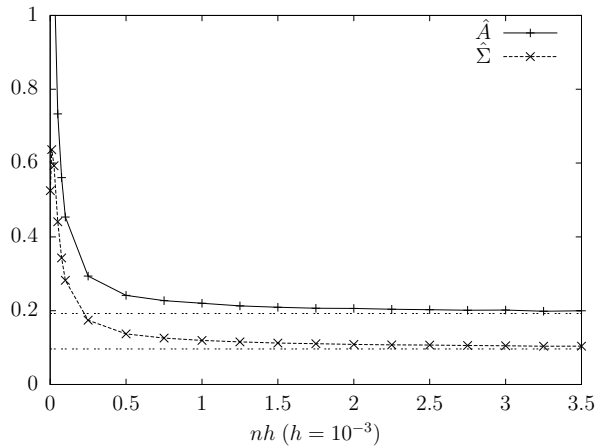}
\caption{Performance of the estimators $\hat{A},\hat{\Sigma}$ in \eqref{eq:sub:brown:effective} as functions
    of the final time $t=nh$ with $h=0.001$.}
  \label{fig:sub:bm:pall}
\end{figure}

Here we study the example
that was originally used in \cite{Pavliotis2007} to illustrate the
failure of classical estimation schemes in the context of multiscale
diffusions for the first time. More precisely we consider the
first-order Langevin equation
\begin{align}
  dx_t = -\nabla V_\alpha\Bigl(x_t,\frac{x_t}{\varepsilon}\Bigr)\,dt + \sqrt{2\sigma}\,dU_t\;,
  \label{eq:brown:grad}
\end{align}
which is a simple model to describe the movement of a Brownian
particle in a two-scale potential $V_\alpha$ subject to thermal
noise -- $U_t$ being a standard Brownian motion. Here we consider
the one-dimensional problem (a two-dimensional example is treated in
see Section \ref{sec:num:multi:langevin}) and further assume that
the two-scale potential is given by a large scale as well as a
fluctuating part: $V_\alpha(x,y) = \alpha V(x) + p(y)$.
Based on these assumptions we can rewrite the Langevin equation as
\begin{align*}
  dx_t = -\Bigl(\alpha V'(x_t) +\frac{1}{\varepsilon}p'\Bigl(\frac{x_t}{\varepsilon}\Bigr)\Bigr)\,dt + \sqrt{2\sigma}\,dU_t\;.
\end{align*}
When the fluctuating part $p$ is sufficiently smooth and periodic with period $L$,
the effective dynamics is given by
\begin{align}
  dX_t = -AV'(X_t)\,dt + \sqrt{2\Sigma}\,dW_t\;,\label{eq:sub:brown:effective}
\end{align}
where the effective coefficients are given by $A=\alpha L^2/(Z_{+}Z_{-})$ and
$\Sigma = \sigma L^2(Z_{+}Z_{-})$, where $Z_{\pm} = \int_0^Le^{\pm p(y)/\sigma}\,dy$, see \cite{Pavliotis2007}
for details. Here we consider $V(x) = x^2/2$ and $p(y) = \cos(y)$, so that \eqref{eq:sub:brown:effective}
is the SDE of an Ornstein-Uhlenbeck process with
\begin{align*}
  A = \frac{\alpha}{I_0(\sigma^{-1})^2}\qquad\text{and}\qquad
  \Sigma = \frac{\sigma}{I_0(\sigma^{-1})^2}\;,
\end{align*}
where $I_{0}(z)$ denotes the modified Bessel function of first kind,
cf.\ \cite[ch.~$9.6$]{Abramowitz1964}. We note that both the
effective drift and the effective diffusion depend on the diffusion
$\sigma$ of the original fast/slow system\footnote{We remark that in
the numerical examples presented in \cite[sec.\
$4.3$]{Frederix2011a} the drift coefficient in the homogenized
equation is assumed to be known when estimating the diffusion
coefficient, even though it depends explicitly on the unknown
$\sigma$.}. Figure \ref{fig:sub:bm:pall} shows the performance of
the estimation scheme when applied to observations of the fast/slow
system with $(\alpha,\sigma) = (1,0.5)$ and $\varepsilon=0.1$.
As for the examples in the previous section both estimators
$\hat{A}$ and $\hat{\Sigma}$ are biased for small final times
$t=nh$. Using longer time series, i.e.\ increasing $nh$, reduces this
bias and both estimators approach the true values (dashed lines)
respectively.

\subsubsection{Brownian motion in a two-scale potential: A quadratic potential in two dimensions}
\label{sec:num:multi:langevin2d}
As a first example to illustrate that the proposed methodology can readily be extended to
multivariate processes, we consider here a generalization of
\eqref{eq:brown:grad} in two dimensions
\begin{align*}
  dx_t = -\nabla V\Bigl(x_t,\frac{x_t}{\varepsilon};M\Bigr)\,dt + \sqrt{2\sigma}\,dU_t\;,
\end{align*}
where $V(\cdot,\cdot;M)$ denotes again a two-scale potential with
$M$ being a set of parameters controlling the drift and $U_t$
denotes a standard two-dimensional Brownian motion. As with the
one-dimensional case, we assume that the two-scale potential
$V(\cdot,\cdot;M)$ is given by a large scale as well as a
fluctuating part, with the fluctuating part being separable:
$V(x,y;M) = V(x;M) + p_1(u) + p_2(v)$, with $x,y\in\R^2$ and
$y=(u,v)^T$. Hence, the original system reads
\begin{align*}
  dx_t = -\Bigl(\nabla V(x_t;M) +
  \frac{1}{\varepsilon}\begin{pmatrix}p_1'\bigl(x_t^1/\varepsilon\bigr)\\p_2'\bigl(x_t^2/\varepsilon\bigr)\end{pmatrix} \Bigr)\,dt + \sqrt{2\sigma}\,dU_t\;,
\end{align*}
with $x_t = (x_t^1,x_t^2)^T\in\R^2$. We take the large scale part to be a quadratic potential
\begin{align*}
  V(x;M) = \frac{1}{2}x^TMx\;,
\end{align*}
with $M$ being symmetric and positive definite, so that the effective dynamics is given by
\begin{align}
  dX_t = -KMX_t\,dt + \sqrt{2\sigma K}\,dW_t\;,\label{eq:sub:brown:effective2d}
\end{align}
for $X_t\in\R^2$ with analytic expressions for $K=\diag(k_1,k_2)$;
see \cite{Pavliotis2007}. With $p_1(u) = \cos(u)$ and $p_2(v) =
\cos(v)/2$ we find
\begin{align*}
  k_1 = \frac{1}{I_0(1/\sigma)^2}\quad\textrm{and}\quad k_2= \frac{1}{I_0(1/(2\sigma))^2}\;,
\end{align*}
where $I_0(z)$ denotes again the modified Bessel function of the
first kind.

\begin{figure}[]
  \centering
  \subfigure[Estimator $\hat{A}$]{
    \includegraphics[width=0.465\textwidth]{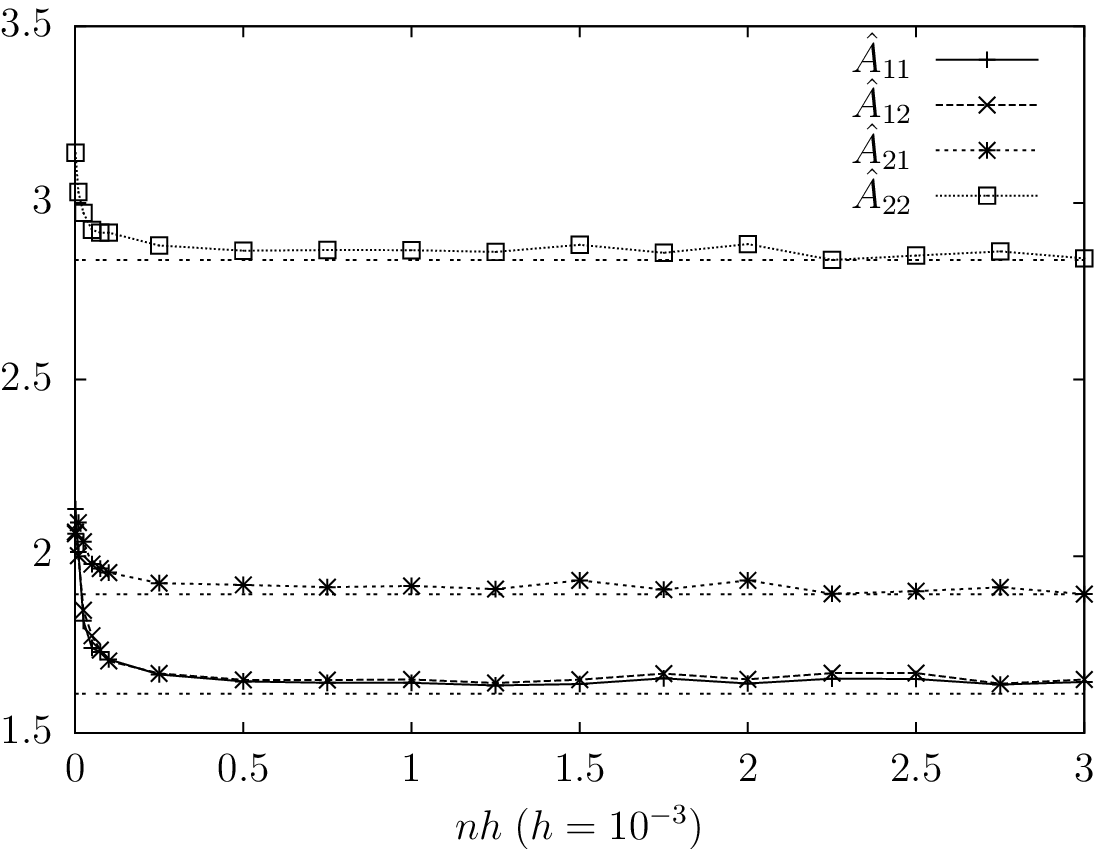}
    \label{fig:sub:bm:p12d}
  }
  \quad
  \subfigure[Estimator $\hat{\Sigma}=\diag(\hat\Sigma_{11},\hat\Sigma_{22})$]{
    \includegraphics[width=0.465\textwidth]{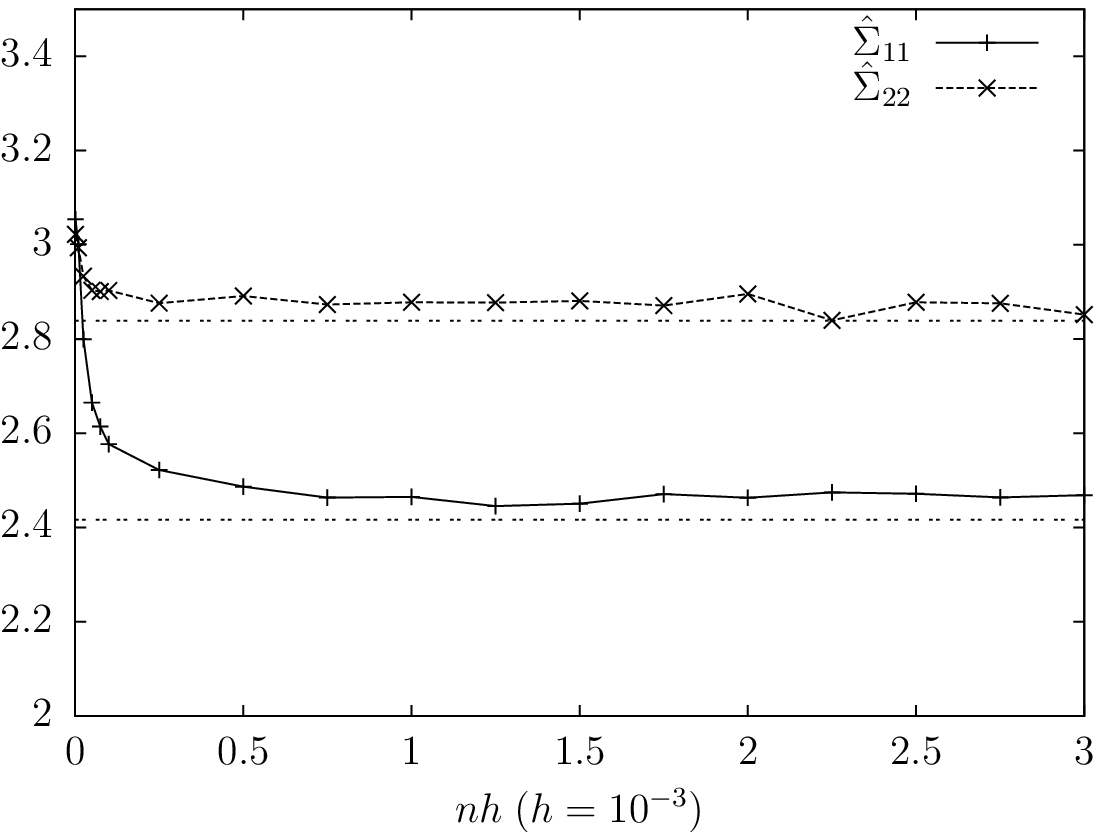}
    \label{fig:sub:bm:p22d}
  }
    \caption{Performance of the estimators $\hat{A},\hat{\Sigma}$ in \eqref{eq:sub:brown:effective2d} as functions
    of the final time $t=nh$ with $h=0.001$.}
  \label{fig:sub:bm:pall2d}
\end{figure}

Since both identities in \eqref{eq:est:start} have pendants for
multivariate processes, the methodology introduced here can be
readily applied to estimate both effective drift matrix $A:=KM$ and
the effective diffusion matrix $\Sigma:=2\sigma K$ in
\eqref{eq:sub:brown:effective2d}. The only difference is that the
system of equations corresponding to \eqref{eq:fun:form:param:sys}
and \eqref{eq:fun:form:param:sys2}, respectively, is a matrix
equation in this case, and similar techniques to obtain the best
approximation (both formally and numerically) can be employed
\cite{Penrose1956}. Figure \ref{fig:sub:bm:pall2d}
depicts the performance of the estimation scheme when applied to
observations of the fast/slow system with $M =
\bigl(\begin{smallmatrix}2&2\\2&3\end{smallmatrix}\bigr)$,
$\sigma=3/2$, and $\varepsilon=0.1$. Figure \ref{fig:sub:bm:p12d}
shows the estimated values of the four effective drift coefficients:
As for the one-dimensional examples, the estimators also show here an
approaching behavior towards the target values when increasing $t=nh$, yielding
an accurate estimate of $A=KM$ for $t=nh$ sufficiently large.
The same behavior when increasing
$t=nh$ is also observed for the estimated (diagonal) effective
diffusion coefficient $\Sigma=2\sigma K$, as it is shown in Figure
\ref{fig:sub:bm:p22d}. We note that although the curves in Figure
\ref{fig:sub:bm:p22d} show a minor gap for larger $t=nh$, the
relative error is less than $2\%$ in both cases for $t\ge 0.75$.

\subsubsection{Eddy diffusivity for the Taylor-Green flow}
\label{sec:num:multi:taylorgreen} 
\begin{figure}[]
    \centering
    \subfigure[time step $h=10^{-3}$]{
      \includegraphics[width=0.465\textwidth]{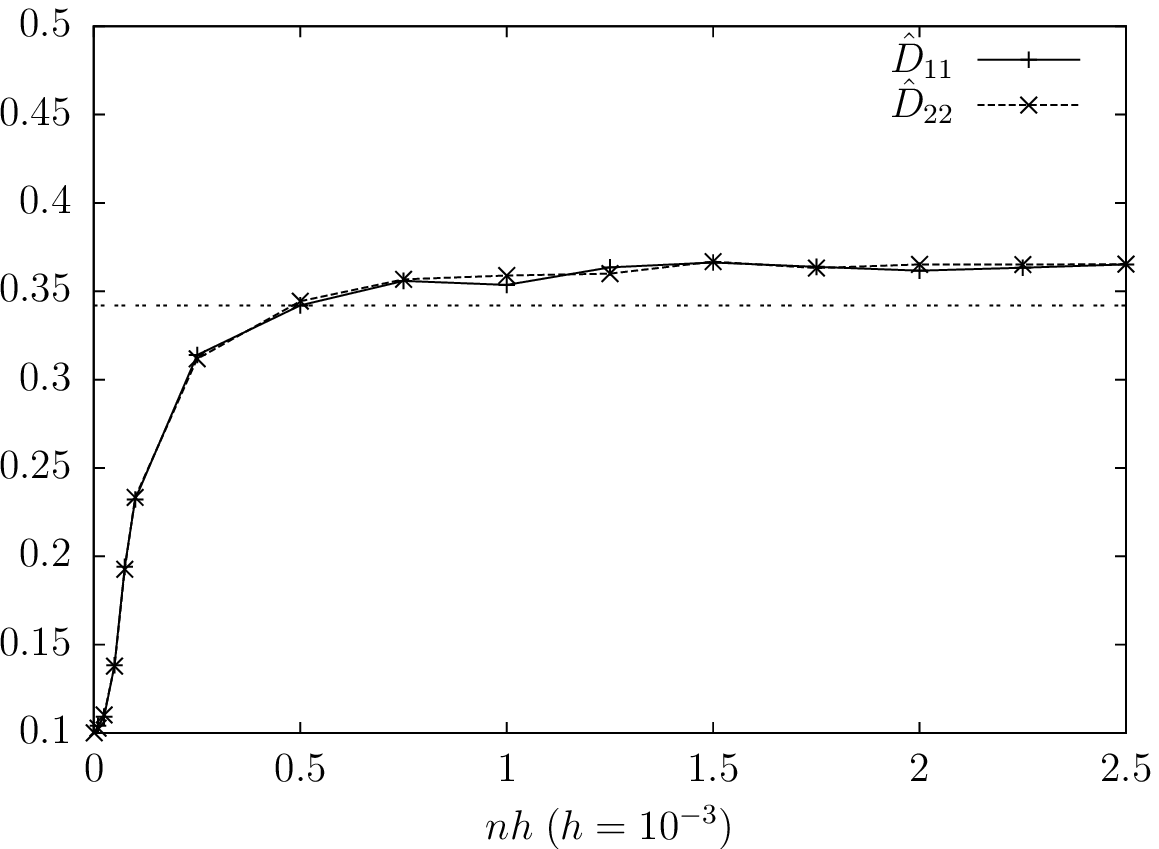}
      \label{fig:sub:tg:normal}
    }
    \quad
    \subfigure[time step $h=10^{-4}$]{
      \includegraphics[width=0.465\textwidth]{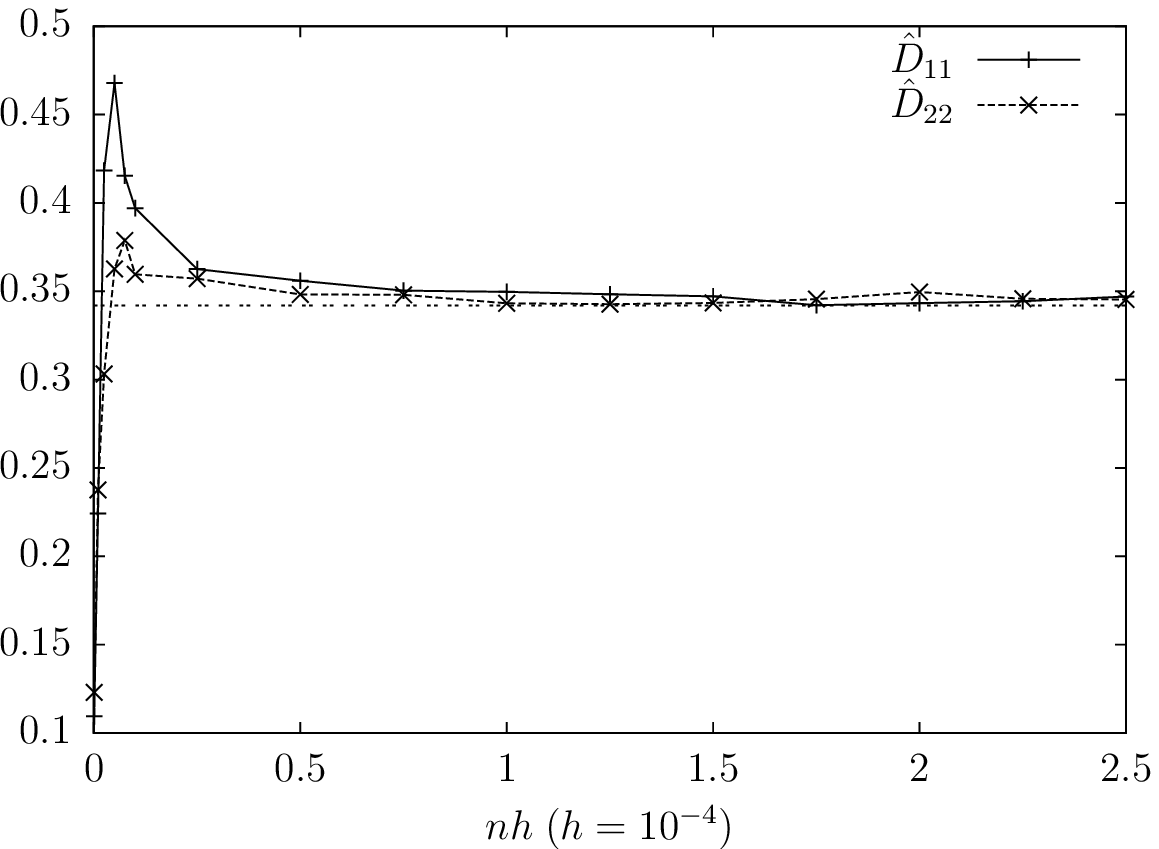}
      \label{fig:sub:tg:fine}
    }
    \caption{Performance of the estimated diagonal elements $\hat{D}_{11}$ and
        $\hat{D}_{22}$ of $\hat{D}$ as functions of the final time $t=nh$.}
    \label{fig:sub:tg:pall}
  \end{figure}

Here we use our estimator to
  estimate the eddy diffusivity (effective diffusion coefficient) of a tracer particle 
  moving in a two-dimensional
  cellular flow and subject to molecular diffusion. This is a very well studied problem and it is
  known that the position of the tracer particle converges, under the diffusive rescaling, to a Brownian
  motion with a diffusion coefficient (covariance matrix) that can be calculated in terms of the
  solution of an appropriate Poisson equation, see  e.g.\ \cite{Majda1999} or \cite[ch.\ $13$]{Pavliotis2008book}.
  The equation for the position of the tracer particle is
  \begin{align*}
    dx_t = v(x_t)\,dt + \sqrt{2\kappa}\,dU_t\;,
  \end{align*}
  where $v$ is a periodic divergence-free velocity field and $\kappa$ denotes the small-scale diffusivity.
  In the numerical simulations below we will take $v$ to
  be the Taylor-Green flow, $v = J \nabla \psi_{TG}$ where $J =\bigl(\begin{smallmatrix}0&-1\\1&0\end{smallmatrix}\bigr)$
  and $\psi_{TG}(u,v) = \sin(u) \sin(v)$. We set $x^{\varepsilon} := \varepsilon x(t/\varepsilon^2)$ to obtain the equation
  \begin{align}
    dx^\varepsilon_t = \frac{1}{\varepsilon}v(x^\varepsilon_t/\varepsilon)\,dt +
    \sqrt{2\kappa}\,dU_t\;.\label{eq:num:tg:rescaled}
  \end{align}
  In the limit as $\varepsilon$ tends to $0$, $x^\varepsilon$ converges weakly to a Brownian motion with diffusion
  tensor $D$, the eddy diffusivity. The goal here is to obtain an estimator $\hat{D}$ of $D$
  using the proposed methodology. It is known that the off-diagonal elements of the eddy diffusivity
  for the Taylor-Green flow vanish, and that the two diagonal elements are equal, and our
  numerical experiments are consistent with these results.
  Figure \ref{fig:sub:tg:pall} shows the performance of diagonal elements of $\hat{D}$
  using \eqref{eq:num:tg:rescaled} with $\varepsilon=0.1$ and $\kappa=0.1$.
  Even though $D$ is not known explicitly as a function of $\kappa$, it can be approximated
  accurately either by solving the Poisson equation using a spectral method or by performing a long time
  Monte Carlo simulation \cite{Cotter2009}. In the aforementioned work the value $\bar{d} = 0.342$ has been
  reported as an approximation of the diagonal elements.
  When a time step $h=10^{-3}$ is used, Figure \ref{fig:sub:tg:normal}
  shows that the estimators behave qualitatively in the same way as in the previous examples: increasing
  $t=nh$ drives the estimators towards limiting values. Noteworthy is that although there
  are differences between these limiting values
  of the estimators and the target value $\bar{d}$ (dashed horizontal line),
  the relative error is less than $7\%$ in both cases for $t=nh\ge 0.5$.
  Moreover, the differences in \ref{fig:sub:tg:normal} are mainly due to temporal discretizations and
  not due to multiscale effects
  as can be verified with Figure \ref{fig:sub:tg:fine} where the performance of the estimators for
  the same experiment but with a smaller time step is shown.

\subsubsection{Truncated Burgers' Equation}
\label{sec:num:multi:burgers}
\begin{figure}[]
  \centering
  \includegraphics[width=0.465\textwidth]{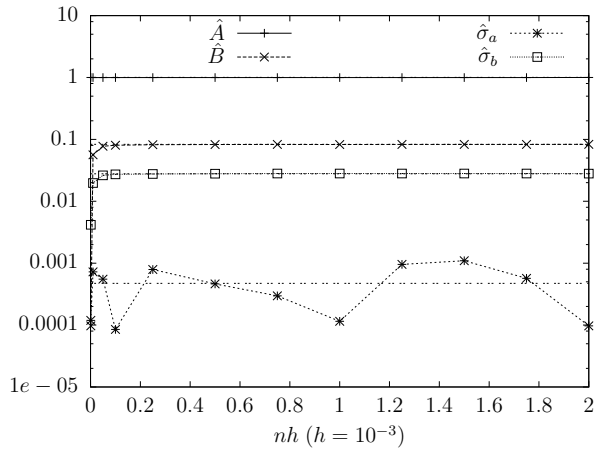}
\caption{Performance of the estimators $\hat{A},\hat{B},\hat{\sigma}_a$, and $\hat{\sigma}_b$ in
  \eqref{eq:num:burgers:effective} as functions of the final time $t=nh$ with $h=0.001$.}
  \label{fig:sub:burgers:pall}
\end{figure}

We consider here an
appropriately rescaled variant of the stochastic Burgers' equation
\begin{align*}
  du_t = \Bigl(\frac{1}{\varepsilon^2}(\partial_x^2 + 1)u_t + \frac{1}{2\varepsilon}\partial_xu_t^2 + \nu u_t \Bigr)\,dt
  + \frac{1}{\varepsilon} \mathcal{Q}\,d\mathcal{W}_t
\end{align*}
on$[0,\pi]$ equipped with homogeneous Dirichlet boundary conditions. Under technical assumptions on the covariance
operator $\mathcal{Q}$ of the space-time white noise $\mathcal{W}$, one can show (see
\cite{Abdulle2012} and references therein for details) that the coefficients of
the three-term truncated representation of the solution have to solve the following multiscale SDE
\begin{align*}
  dx_t &= \bigl(\nu x_t - \frac{1}{2\varepsilon}(x_ty_t^1+y_t^1y_t^2)\bigr)\,dt\;,\\
  dy_t^1 &= \bigl(\nu y_t^1-\frac{3}{\varepsilon^2}y_t^1- \frac{1}{2\varepsilon}(2x_ty_t^2-{x_t}^2)\bigr)\,dt +
  \frac{q_1}{\varepsilon}\,dV_t^1\;,\\
  dy_t^2 &= \bigl(\nu y_t^2 - \frac{8}{\varepsilon^2}y_t^2 + \frac{3}{2\varepsilon}x_ty_t^1\bigr)\,dt +
  \frac{q_2}{\varepsilon}\,dV_t^2\;,
\end{align*}
with $V_t^1$ and $V_t^2$ being independent standard Brownian motions.
Therein the covariance operator $\mathcal{Q}$ is such
that noise acts only on the fast modes directly.
For the truncated system standard homogenization theory applies and yields
\begin{align}
dX_t = \bigl(AX_t-B{X_t}^3\bigr)\,dt + \sqrt{\sigma_a+\sigma_b{X_t}^2}\,dW_t\label{eq:num:burgers:effective}
\end{align}
as the effective dynamics with true parameters
\begin{align*}
A = \nu + \frac{{q_1}^2}{396}+\frac{{q_2}^2}{352}\;,\quad B = \frac{1}{12}\;,\quad
\sigma_a = \frac{{q_1}^2{q_2}^2}{2112}\;,\text{ and}\quad \sigma_b = \frac{{q_1}^2}{36}\;.
\end{align*}
See \cite{Blomker2007} for details.
Figure \ref{fig:sub:burgers:pall} shows the performance of the estimation scheme when applied to
observations of the three dimensional fast/slow system with $\nu = 1$, $(q_1,q_2) = (1,1)$, and
$\varepsilon=0.1$.
Since the true values of the effective coefficients (dashed lines)
are of different orders for these particular choices, a
semi-logarithmic scale is adopted for the sake of clarity. The plots
show qualitatively the same behavior as in the previous examples
when increasing the final time $t=nh$ and suggest that the
estimation procedure yields accurate estimators. Only the estimated
value $\hat{\sigma}_a$ fluctuates around the true value. We
note that the true value $\sigma_a$ is very small ($\approx 5\cdot
10^{-4}$) so that this coefficient has only marginal influence in
the complete diffusion function. Furthermore, recall that even the
time-step in the Euler-Maruyama discretization is larger
($h=10^{-3}$), so that the fluctuations are expected to be due to
discretization errors. As a matter of fact, considering a finer
discretization (i.e.\ increasing $m,n$, and $N$) reduces the
fluctuations (not shown here).

\subsubsection{Fast chaotic noise}
\label{sec:num:multi:chaos1}
Our methodology is also applicable to a system of ordinary differential equations (ODEs)
where the stochastic noise is replaced by deterministic chaos. In particular, we will
consider an ODE driven by one of the components of an appropriately rescaled Lorentz system.
More precisely, consider as an example the following system
\begin{subequations}
\begin{align}
  \frac{dx}{dt} & = x-x^3 + \frac{\lambda}{\varepsilon}(1+\nu x^2) y_2\;,\label{eq:num:chaos:fastslow:slow}\\
  \frac{dy_1}{dt} & = \frac{10}{\varepsilon^2}(y_2-y_1)\;,\label{eq:num:chaos:fastslow:fast1}\\
  \frac{dy_2}{dt} & = \frac{1}{\varepsilon^2}(28y_1 - y_2 - y_1y_3)\;,\label{eq:num:chaos:fastslow:fast2}\\
  \frac{dy_3}{dt} & = \frac{1}{\varepsilon^2}(y_1y_2 - \frac{8}{3}y_3)\;,\label{eq:num:chaos:fastslow:fast3}
\end{align}
\label{eq:num:chaos:fastslow}%
\end{subequations}
where the fast component $y = (y_1,y_2,y_3)^T$ solves the Lorenz equation. In the sequel we investigate two
different couplings between the fast and the slow process by choosing $\nu\in\{0,1\}$.

According to \cite[ch.~$11.7.2$]{Pavliotis2008book} (see also \cite[ex.~$6.2$]{Givon2004}),
when eliminating the fast chaotic variable $y$, the approximate dynamics for $\nu=0$ is given by
\begin{align}
  dX_t = A\bigl(X_t-{X_t}^3\bigr)\,dt + \sqrt{\sigma}\,dW_t\;,\label{eq:num:chaos:effective}
\end{align}
with $A=1$ and the diffusion coefficient that is given by the Green-Kubo formula
\begin{align}
  \sigma = 2\lambda^2\int_0^\infty\lim_{T\rightarrow\infty}\frac{1}{T}\int_0^T\psi^s(y)\psi^{s+t}(y)\,ds\,dt\;.
  \label{eq:num:chaos:sigma}
\end{align}
Therein $\psi^t(y) = e_2\cdot\varphi^t(y)$ with $\varphi^t(y)$
denoting the solution of the fast process $y$ at time $t$ when
$\varepsilon=1$ and $e_2 = (0,1,0)^T$. The convergence of the
solution of \eqref{eq:num:chaos:fastslow:slow} to the solution of
\eqref{eq:num:chaos:effective} can be justified rigorously using the
recent results from \cite{Melbourne2011}. However, the above
expression for $\sigma$ is not useful practically: Not only does it not
give an analytical value for $\sigma$, also using it to obtain
$\sigma$ numerically is computationally expensive.
Hence, it would
be advantageous to use the methodology proposed here and estimate
the effective coefficients via observations of the complete
(deterministic) fast/slow system. To illustrate numerically that our
estimation procedure can indeed deal with this problem, we apply it
to observations of the deterministic fast/slow system using $\lambda
= 2/45$ and $\nu=0$ to estimate both the drift
and the diffusion coefficients. Since the
system is deterministic, classical solvers for ODEs may be employed.
For example, depending on the stiffness (i.e.\ on $\varepsilon$) of
the system, either a fourth order Runge-Kutta scheme or a solver
based on numerical differentiation formulas is used.
\begin{figure}[]
  \centering
  \subfigure[Estimator $\hat{A}$]{
    \includegraphics[width=0.465\textwidth]{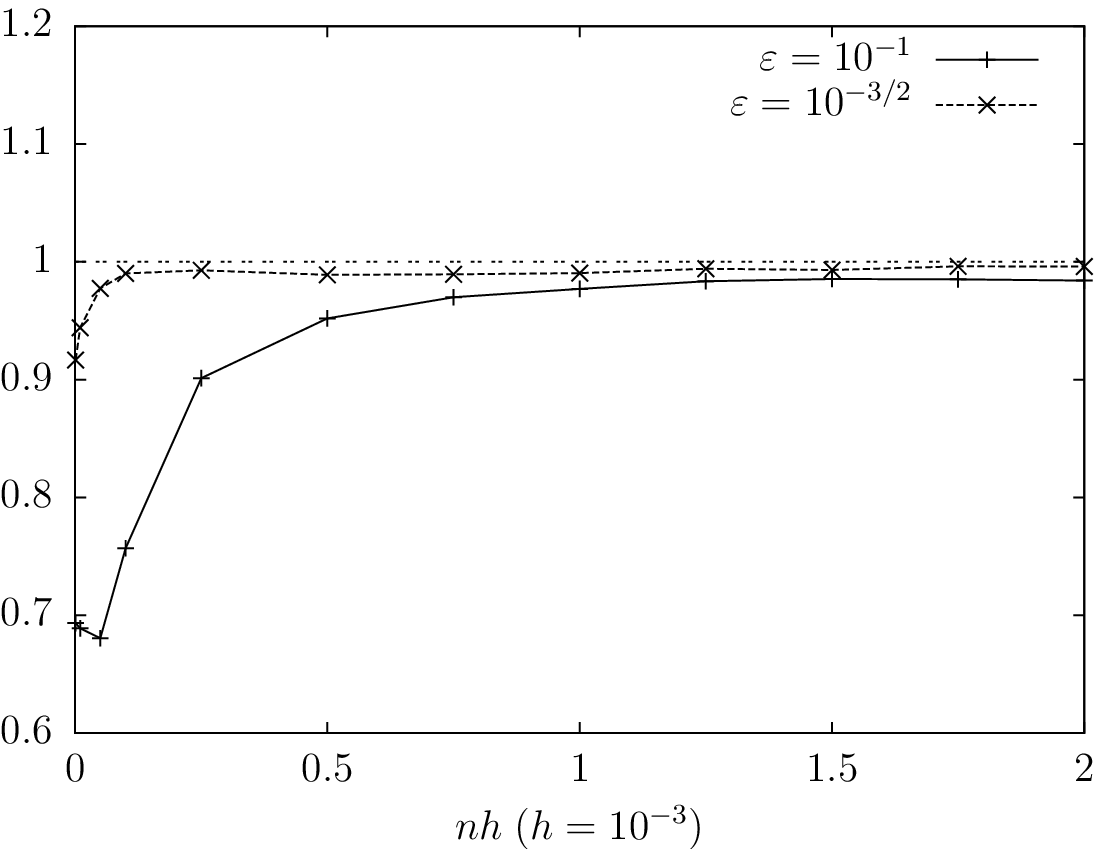}
    \label{fig:sub:chaos:p1}
  }
  \quad
  \subfigure[Estimator $\hat{\sigma}$]{
    \includegraphics[width=0.465\textwidth]{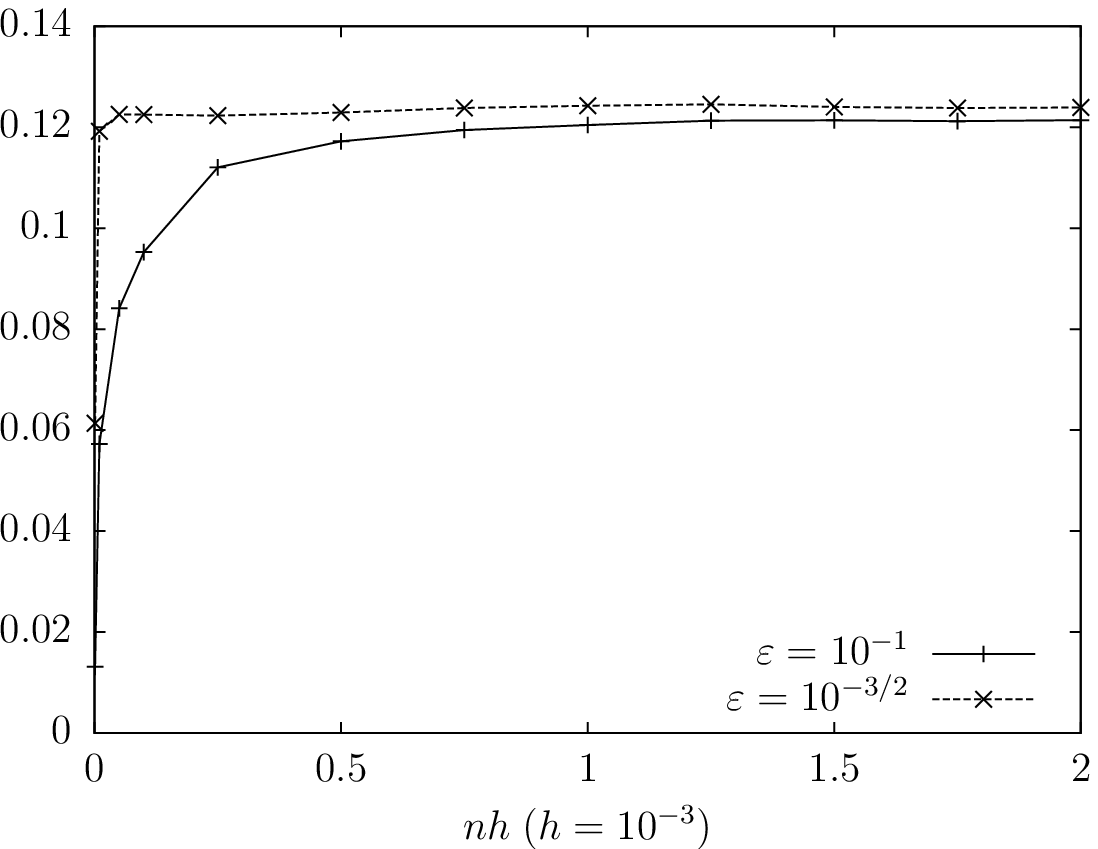}
    \label{fig:sub:chaos:p2}
  }
  \caption[]{Performance of the estimation scheme applied to the deterministic system
    \eqref{eq:num:chaos:fastslow} using $\lambda = 2/45$, $\nu=0$, and $\varepsilon\in\{10^{-3/2},10^{-1}\}$.
    The final time of the considered time series is $t=nh$ with $h=0.001$ and the true effective drift parameter
    is $A = 1$.}
  \label{fig:sub:chaos}
\end{figure}
Figure \ref{fig:sub:chaos} depicts the estimated values as functions
of the final time $t=nh$ for two different choices of the scale
separation $\varepsilon\in\{10^{-3/2},10^{-1}\}$.
The value $\varepsilon=10^{-3/2}$ is the same with that used
in \cite{Givon2004}, thus allowing for direct comparisons to be
made; we will return to this point shortly. The estimated drift
parameter $\hat{A}$ (Figure \ref{fig:sub:chaos:p1}) shows (for both
values of $\varepsilon$) the typical behavior we expect in the
context of multiscale diffusions. While the estimator is biased for
small values of $nh$, increasing $nh$ significantly reduces the bias
so that the estimator approaches the true value (dashed line). The
estimators of the effective diffusion coefficient $\hat{\sigma}$
(Figure \ref{fig:sub:chaos:p2}) also approach a limiting value when
increasing $nh$, with minor fluctuations for both values of
$\varepsilon$. In both plots, one observes a performance difference
of the estimators for different values of $\varepsilon$. In fact,
the more distinctive the scale separation (i.e.\ the smaller
$\varepsilon$) between fast and slow components, the faster the
estimators approach a limiting value.

For both drift and diffusion estimator
 the curves for different values of
$\varepsilon$ give slightly different values
\begin{align*}
  \hat{A}_\varepsilon\approx\begin{cases} 0.984 &\textrm{, if } \varepsilon=10^{-1} \\
    0.998 & \textrm{, if } \varepsilon=10^{-3/2} \end{cases}
  \qquad\textrm{and}\qquad
  \hat{\sigma}_\varepsilon \approx\begin{cases} 0.121 &\textrm{, if } \varepsilon=10^{-1} \\
    0.124 & \textrm{, if } \varepsilon=10^{-3/2} \end{cases}
\end{align*}
at time $nh=2$. Hence for $\varepsilon
>0$ there exists an additional bias, which is the reason for the
observed difference in the estimators for different $\varepsilon$.
On the other hand, one expects that as $\varepsilon$ decreases the
estimators of the effective coefficients become more accurate, as
one observes here for the estimated drift $\hat{A}$. Thus,
$\hat\sigma$ is also expected to be more accurate as $\varepsilon$
decreases. Since no analytic formula for the diffusion coefficient
exists, the estimator $\hat{\sigma}$ is compared with
alternative numerical approximations available in \cite[ex.~$6.2$
and ill.~$10.5$]{Givon2004}. In the numerical experiments performed
in this study the value $\varepsilon = 10^{-3/2}$ was adopted
giving a value of $0.126\pm 0.003$ using Gaussian (second) moment
approximations based on a modified Euler-Maruyama discretization of
the effective dynamics, and a value of $0.13\pm 0.01$ based on the
heterogeneous multiscale method (HMM) with a discretization of the
Green-Kubo formula for the effective diffusion coefficient; see
also \cite{Fatkullin2004} for more elaborated HMM based numerical
schemes applied to the Lorenz 96 model. Thus, we have a very good
agreement of the result obtained by the estimation proposed here
with these previously reported values. It is worth to mention that,
unlike the procedure introduced here, the methods employed to
determine the effective diffusion coefficient in \cite{Givon2004}
assume that the effective drift parameter $A$ is known. While the
HMM can easily be adapted to the case of an unknown drift parameter,
it is not straightforward to incorporate the unknown drift parameter
in the estimation based on Gaussian moment approximations. In any
event, incorporating the drift estimation would yield an even larger
statistical error for these methods, while the results based on the
presented scheme show only minor fluctuations; see Figure
\ref{fig:sub:chaos}.

\begin{figure}[]
  \centering
  \subfigure[Drift Parameters]{
    \includegraphics[width=0.465\textwidth]{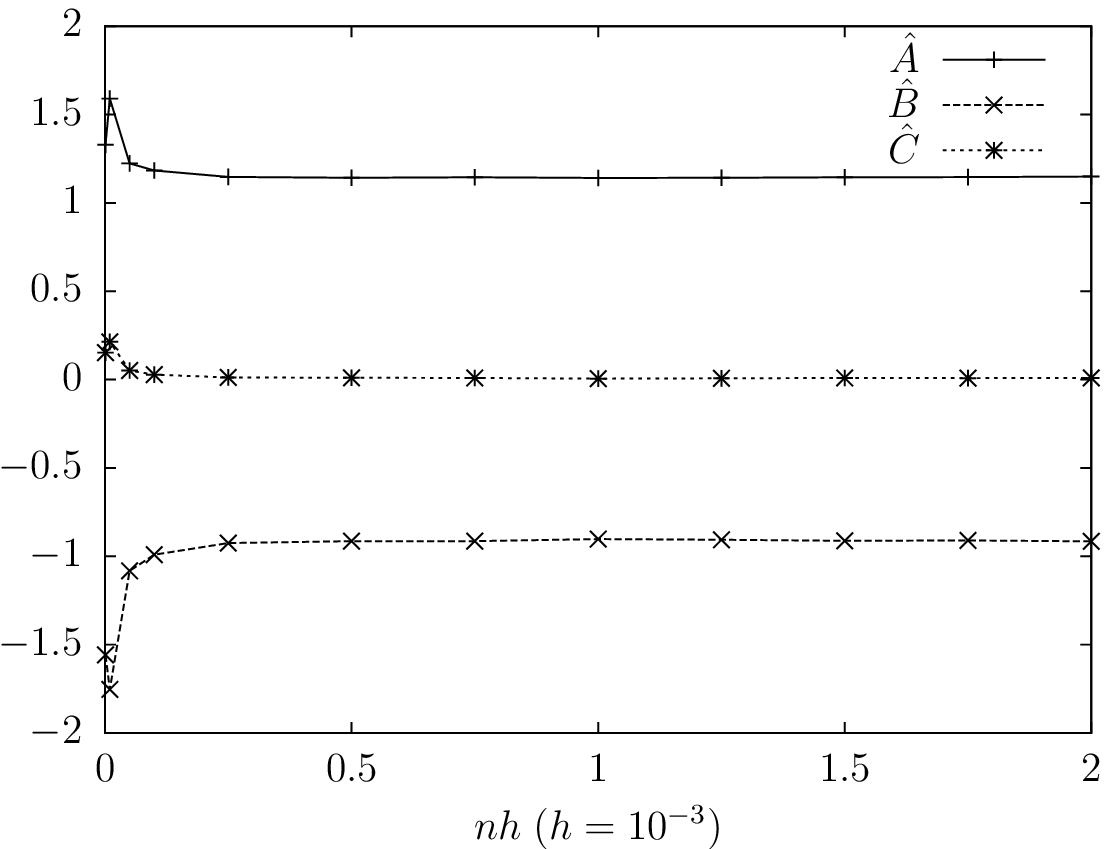}
    \label{fig:sub:chaos:mod:drift}
  }
  \quad
  \subfigure[Diffusion Parameters]{
    \includegraphics[width=0.465\textwidth]{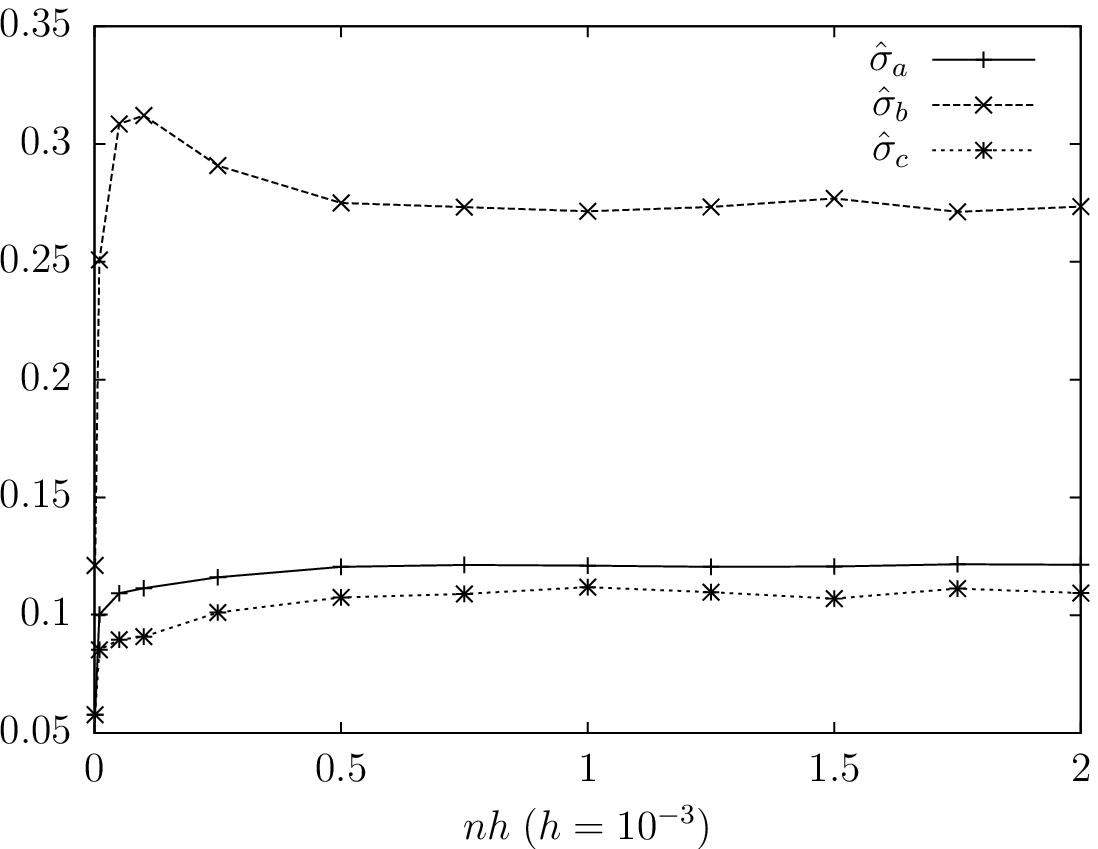}
    \label{fig:sub:chaos:mod:diffusion}
  }
  \caption[]{Performance of the estimation scheme applied to the deterministic system
    \eqref{eq:num:chaos:fastslow} using $\lambda = 2/45$, $\nu=1$, and $\varepsilon=10^{-3/2}$. The final time of
    the considered time series is $t=nh$ with $h=0.001$.}
  \label{fig:sub:chaos:mod}
\end{figure}

Choosing $\nu=1$ in \eqref{eq:num:chaos:fastslow:slow} and following
the methodology outlined in \cite[ch.~$11$]{Pavliotis2008book},
yields the following effective dynamics
\begin{align}
  dX_t = \bigl(AX_t+B{X_t}^3+C{X_t}^5\bigr)\,dt + \sqrt{\sigma_a+\sigma_b {X_t}^2+\sigma_c{X_t}^4}\,dW_t\;.
  \label{eq:num:chaos:effective2}
\end{align}
The effective coefficients are now given by
\begin{align*}
  A = 1+\sigma\;,\quad B=\sigma-1\;,\quad C=0\;,\footnotemark[6]
  \quad \sigma_a = \sigma\;,\quad \sigma_b = 2\sigma\;,\quad
  \sigma_c = \sigma\;,
\end{align*}
\footnotetext[6]{Although we know theoretically that $C=0$, we estimate $C$ nonetheless with the
  proposed scheme to illustrate that the novel scheme can correctly identify the relevant parameters
  in a model that contains more parameters than necessary, for instance given by a Taylor or Fourier
  series expansion.}
with $\sigma$ being as in \eqref{eq:num:chaos:sigma}. We apply our
methodology again to observations of the deterministic fast/slow
system using $\lambda = 2/45$, $\nu=1$, and $\varepsilon=10^{-3/2}$
to estimate all six effective coefficients. Figure
\ref{fig:sub:chaos:mod} illustrates the estimated values of both
drift and diffusion parameters as functions of the final time
$t=nh$.
We observe the procedure's typical behavior in the context of
multiscale observations. In fact, for both drift
(Figure \ref{fig:sub:chaos:mod:drift}) and diffusion parameters
(Figure \ref{fig:sub:chaos:mod:diffusion}) the estimators approach limiting
values with only minor fluctuations when increasing $nh$.


%% file: outlook.tex
\section{Conclusion}
\label{sec:outlook}

We have developed a numerical methodology for estimating multiple
parameters in a coarse-grained equation (in one or multiple
dimensions) based on observation from an associated fast/slow system
that posses a multiscale structure. This problem is far from
straightforward, not only due to the multiscale effects present in
the available data, but also due to difficulties associated with
estimating parameters when both the drift and the diffusion coefficients
are state dependent.

The approach developed in this study combines a number of different techniques.
On the one hand, the derivation of the estimators relies on simple
identities based on the martingale property for stochastic integrals
and It{\^o} isometry. On the other, we exploit our freedom in varying
the initial condition in combination with standard techniques from
inverse problems to define the parameter estimators via best
approximation.

We demonstrated via a detailed numerical study that the proposed
inference scheme provides us with accurate estimates for parameters
in both the drift and the diffusion coefficients in systems with multiscale
structure and state dependent noise. In fact, the proposed methodology
appears to be accurate and effective even when the stochastic noise in
the system is replaced by deterministic chaos.

While the feasibility study of the parameter estimation for
multiscale diffusions and the initial presentation of the estimation
scheme is the main development here, clearly many open problems and
questions remain to be addressed. One such open question is the
rigorous analysis of the algorithm to investigate its asymptotic
properties and to scrutinize its limitations. Furthermore, the
rigorous analysis of the algorithm is expected to reveal insights
into the dependency on the algorithm-defining parameters that in
turn can be used to reduce the computational complexity of the
methodology.

Clearly, from a practical point of view there are different
strategies to improve the efficiency of the estimators. A first
starting-point could be the usage of techniques with an accelerated
convergence instead of the brute-force Monte Carlo sampling to
approximate the involved expectations, e.g.\ quasi Monte Carlo
\cite{Niederreiter1992} or variance reduction techniques \cite[ch.\
$16$]{Kloeden1992}. Also recent work on Multilevel Monte Carlo
methods \cite{Giles2008} appears very appealing in this context,
although care might have to be taken due to the nonlinear nature of
the underlying model SDE; cf.\ \cite{Hutzenthaler2012pre}.

Several questions also arise naturally within the presented
framework of varying the initial condition to set up a system of
equations: How many initial conditions need to be considered? Where
to locate the initial conditions: Equispaced or distributed
differently? How does the choice of the initial condition influence
the accuracy of the estimator? In fact, preliminary numerical
experiments suggest that an alternative distribution of the initial
conditions improves the accuracy.
Hence, an ``optimal'' distribution of the initial condition
is expected to reduce the computational cost considerably.

Another interesting point concerns the best approximation.
Here, we defined the estimator as the element that minimizes the
residual of a linear system of equations with respect to the
Euclidean norm. Thus, natural questions are, e.g.\ as to whether an
alternative norm might be more appropriate, and how regularizing the
minimization problem (e.g.\ by a truncated singular value
decomposition or a Tikhonov regularization) affects the estimator.
Also methodologies that ``optimize'' the linear system of equations
to obtain more accurate best approximations, as they are for
instance used in the reconstruction of tomographic problems
\cite{Lustfeld2011}, are appealing.

From a computational point of view, employing parallel computing
strategies seem also promising. In fact, the presented algorithm
consists of multiple parts (mainly when generating the observations)
without mutual dependency.
Thus, these parts lead to problems that are straightforward to parallelize.

There are applications, however, where only one long time-series is
  available, rather than several short ones, as required for the 
  methodology developed here to estimate the drift and diffusion coefficients. 
  It is an interesting question whether this methodology
  can be extended to cases where a long time-series is available only. This and
  related issues will be treated in feature studies.


%% file: KPK-MMS2011.bbl
\begin{thebibliography}{10}

\bibitem{Abdulle2012}
{\sc A.~Abdulle and G.~A. Pavliotis}, {\em {N}umerical {M}ethods for
  {S}tochastic {P}artial {D}ifferential {E}quations with {M}ultiple {S}cales},
  J. Comput. Phys., 231 (2012), pp.~2482--2497.

\bibitem{Abramowitz1964}
{\sc M.~Abramowitz and I.~A. Stegun}, {\em {H}andbook of {M}athematical
  {F}unctions with {F}ormulas, {G}raphs and {M}athematical {T}ables}, Dove, New
  York, 1964.

\bibitem{Azencott2010}
{\sc R.~Azencott, A.~Beri, and I.~Timofeyev}, {\em {A}daptive {S}ub-sampling
  for {P}arametric {E}stimation of {G}aussian {D}iffusions}, J. Stat. Phys.,
  139 (2010), pp.~1066--1089.

\bibitem{Azencott2011}
\leavevmode\vrule height 2pt depth -1.6pt width 23pt, {\em {P}arametric
  {E}stimation of {S}tationary {S}tochastic {P}rocesses under {I}ndirect
  {O}bservability}, J. Stat. Phys., 144 (2011), pp.~150--170.

\bibitem{Ben-Israel2003}
{\sc A.~Ben-Israel and T.~N.~E. Greville}, {\em {G}eneralized {I}nverses:
  {T}heory and {A}pplications}, CMS Books in Mathematics, Springer, 2003.

\bibitem{Blomker2007}
{\sc D.~Bl{\"o}mker, M.~Hairer, and G.~A. Pavliotis}, {\em {M}ultiscale
  {A}nalysis for {S}tochastic {P}artial {D}ifferential {E}quations with
  {Q}uadratic {N}onlinearities}, Nonlinearity, 20 (2007), pp.~1721--1744.

\bibitem{Calderon2007}
{\sc C.~P. Calderon}, {\em {F}itting {E}ffective {D}iffusion {M}odels to {D}ata
  {A}ssociated with a {"Glassy"} {P}otential: {E}stimation, {C}lassical
  {I}nference {P}rocedures, and {S}ome {H}euristics}, Multiscale Model. Simul.,
  6 (2007), pp.~656--687.

\bibitem{Chauviere2010}
{\sc A.~Chauvi{\`e}re, L.~Preziosi, and C.~Verdier}, eds., {\em {C}ell
  {M}echanics: {F}rom {S}ingle {S}cale-{B}ased {M}odels to {M}ultiscale
  {M}odeling}, Mathematical \& Computational Biology Series, Chapman \&
  Hall/CRC, 2010.

\bibitem{Cotter2009}
{\sc C.~J. Cotter and G.~A. Pavliotis}, {\em Estimating eddy diffusivities from
  noisy {L}agrangian observations}, Commun. Math. Sci., 7 (2009), pp.~805--838.

\bibitem{Cross1993}
{\sc M.~C. Cross and P.~C. Hohenberg}, {\em Pattern formation outside of
  equilibrium}, Rev. Modern Phys., 65 (1993), pp.~851--1112.

\bibitem{Cruz-Uribe2002}
{\sc D.~Cruz-Uribe and C.~J. Neugebauer}, {\em {S}harp {E}rror {B}ounds for the
  {T}rapezoidal {R}ule and {S}impson's {R}ule}, JIPAM. J. Inequal. Pure Appl.
  Math., 3 (2002), pp.~1--22.

\bibitem{E2005}
{\sc W.~E, D.~Liu, and E.~Vanden-Eijnden}, {\em Analysis of multiscale methods
  for stochastic differential equations}, Comm. Pure Appl. Math., 58 (2005),
  pp.~1544--1585.

\bibitem{Einstein1956}
{\sc A.~Einstein}, {\em {I}nvestigations on the {T}heory of {B}rownian
  {M}ovement}, Dover Publications, 1956.
\newblock edited with notes by R. F{\"u}rth, translated by A. D. Cowper.

\bibitem{Fatkullin2004}
{\sc I.~Fatkullin and E.~Vanden-Eijnden}, {\em A computational strategy for
  multiscale systems with applications to lorenz 96 model}, J. Comput. Phys.,
  200 (2004), pp.~605--638.

\bibitem{Fish2009}
{\sc J.~Fish}, {\em {M}ultiscale {M}ethods: {B}ridging the {S}cales in
  {S}cience and {E}ngineering}, Oxford University Press, 2009.

\bibitem{Frederix2011a}
{\sc Y.~Frederix and D.~Roose}, {\em {D}rift-{F}iltered {A}pproach to
  {D}iffusion {E}stimation for {M}ultiscale {P}rocesses}, in Coping with
  Complexity: Model Reduction and Data Analysis, A.~N. Gorban and D.~Roose,
  eds., vol.~75 of Lecture Notes in Computational Science and Engineering,
  Springer, 2011, pp.~269--286.

\bibitem{Frederix2011}
{\sc Y.~Frederix, G.~Samaey, and D.~Roose}, {\em {A}n analysis of noise
  propagation in the multiscale simulation of coarse {F}okker-{P}lanck
  equations}, ESAIM Math. Model. Numer. Anal., 45 (2011), pp.~541--561.

\bibitem{Giles2008}
{\sc M.~B. Giles}, {\em {M}ultilevel {M}onte {C}arlo path simulation}, Oper.
  Res., 56 (2008), pp.~607--617.

\bibitem{Givon2004}
{\sc D.~Givon, R.~Kupferman, and A.~M. Stuart}, {\em Extracting macroscopic
  dynamics: model problems and algorithms}, Nonlinearity, 17 (2004),
  pp.~55--127.

\bibitem{Golub1996}
{\sc G.~H. Golub and C.~F. van Loan}, {\em Matrix Computations}, Johns Hopkins
  Studies in the Mathematical Sciences, Johns Hopkins University Press,
  3rd~ed., 1996.

\bibitem{Griebel2007}
{\sc M.~Griebel, S.~Knapek, and G.~W. Zumbusch}, {\em {N}umerical {S}imulation
  in {M}olecular {D}ynamics: {N}umerics, {A}lgorithms, {P}arallelization,
  {A}pplications}, Texts in Computational Science and Engineering, Springer,
  2007.

\bibitem{Hamilton1994}
{\sc J.~D. Hamilton}, {\em {T}ime {S}eries {A}nalysis}, Princeton University
  Press, 1994.

\bibitem{Hansen1982}
{\sc L.~P. Hansen}, {\em {L}arge {S}ample {P}roperties of {G}eneralized
  {M}ethod of {M}oments {E}stimators}, Econometrica, 50 (1982), pp.~1029--1054.

\bibitem{Horstemeyer2010}
{\sc M.~F. Horstemeyer}, {\em {M}ultiscale {M}odeling: {A} {R}eview}, in
  Practical Aspects of Computational Chemistry, J.~Leszczynski and M.~K.
  Shukla, eds., Springer, 2010, pp.~87--135.

\bibitem{Horsthemke1984}
{\sc W.~Horsthemke and R.~Lefever}, {\em {N}oise-{I}nduced {T}ransitions:
  {T}heory and {A}pplications in {P}hysics, {C}hemistry, and {B}iology},
  vol.~15 of Springer Series in Synergetics, Springer, 1984.

\bibitem{Huerre1998}
{\sc P.~Huerre and M.Rossi}, {\em Hydrodynamic instabilities in open flows}, in
  {H}ydrodynamic and {N}onlinear {I}nstabilities, C.~Godr\`eche and
  P.~Manneville, eds., Cambridge University Press, 1998, pp.~81--294.

\bibitem{Hutzenthaler2012pre}
{\sc M.~Hutzenthaler, A.~Jentzen, and P.~E. Kloeden}, {\em {D}ivergence of the
  multilevel {M}onte {C}arlo {E}uler method for nonlinear stochastic
  differential equations}.
\newblock To appear in Ann. Appl. Probab. (2013).

\bibitem{Kalliadasis2000}
{\sc S.~Kalliadasis}, {\em {N}onlinear instability of a contact line driven by
  gravity}, J. Fluid Mech., 413 (2000), pp.~355--378.

\bibitem{Kalliadasis2012}
{\sc S.~Kalliadasis, C.~Ruyer-Quil, B.~Scheid, and M.~G. Velarde}, {\em
  {F}alling {L}iquid {F}ilms}, vol.~176 of Applied Mathematical Sciences,
  Springer, 2012.

\bibitem{Kevrekidis2004}
{\sc I.~G. Kevrekidis, C.~W. Gear, and G.~Hummer}, {\em {E}quation-free: {T}he
  computer-aided analysis of complex multiscale systems}, AIChE J., 50 (2004),
  pp.~1346--1355.

\bibitem{Kevrekidis2003}
{\sc I.~G. Kevrekidis, C.~W. Gear, J.~M. Hyman, P.~G. Kevrekidis, O.~Runborg,
  and C.~Theodoropoulos}, {\em {E}quation-{F}ree, {C}oarse-{G}rained
  {M}ultiscale {C}omputation: {E}nabling {M}icroscopic {S}imulators to
  {P}erform {S}ystem-{L}evel {A}nalysis}, Commun. Math. Sci., 1 (2003),
  pp.~715--762.

\bibitem{Kevrekidis2009}
{\sc I.~G. Kevrekidis and G.~Samaey}, {\em {E}quation-{F}ree {M}ultiscale
  {C}omputation: {A}lgorithms and {A}pplications}, Annu. Rev. Phys. Chem., 60
  (2009), pp.~321--344.

\bibitem{Kloeden1992}
{\sc P.~E. Kloeden and E.~Platen}, {\em {N}umerical {S}olution of {S}tochastic
  {D}ifferential {E}quations}, Springer, 1992.

\bibitem{Krylov1999}
{\sc N.~V. Krylov}, {\em {O}n {K}olmogorov's equations for finite dimensional
  diffusions}, in Stochastic {PDE}'s and {K}olmogorov equations in infinite
  dimensions ({C}etraro, 1998), vol.~1715 of Lecture Notes in Mathematics,
  Springer, 1999, pp.~1--63.

\bibitem{Kuramoto2003}
{\sc Y.~Kuramoto}, {\em {C}hemical {O}scillations, {W}aves, and {T}urbulence},
  Chemistry Series, Dover Publications, 2003.
\newblock Slightly corrected republication of the work originally published by
  Springer in 1984.

\bibitem{Kutoyants2004}
{\sc Y.~A. Kutoyants}, {\em {S}tatistical {I}nference for {E}rgodic {D}iffusion
  {P}rocesses}, Springer, 2004.

\bibitem{Liptser2010}
{\sc R.~S. Liptser and A.~N. Shiryaev}, {\em {S}tatistics of {R}andom
  {P}rocesses: {I}. {G}eneral {T}heory}, Stochastic Modelling and Applied
  Probability Series, Springer, 2nd~ed., 2010.
\newblock Translated by A. B. Aries.

\bibitem{Lo1988}
{\sc A.~W. Lo}, {\em {M}aximum {L}ikelihood {E}stimation of {G}eneralized
  {I}t\^o {P}rocesses with {D}iscretely {S}ampled {D}ata}, Econometric Theory,
  4 (1988), pp.~231--247.

\bibitem{Lustfeld2011}
{\sc H.~Lustfeld, J.~A. Hirschfeld, M.~Rei{\ss}el, and B.~Steffen}, {\em
  Enhancement of precision and reduction of measuring points in tomographic
  reconstructions}, Phys. Lett. A, 375 (2011), pp.~1167--1171.

\bibitem{Majda2008}
{\sc A.~J. Majda, C.~Franzke, and B.~Khouider}, {\em An applied mathematics
  perspective on stochastic modelling for climate}, Philos. Trans. R. Soc.
  Lond. Ser. A Math. Phys. Eng. Sci., 366 (2008), pp.~2427--2453.

\bibitem{Majda1999}
{\sc A.~J. Majda and P.~R. Kramer}, {\em Simplified models for turbulent
  diffusion: theory, numerical modelling, and physical phenomena}, Phys. Rep.,
  314 (1999), pp.~237--574.

\bibitem{Mazo2002}
{\sc R.~M. Mazo}, {\em {B}rownian {M}otion: {F}luctuations, {D}ynamics, and
  {A}pplications}, International Series of Monographs on Physics, Oxford
  University Press, 2002.

\bibitem{Melbourne2011}
{\sc I.~Melbourne and A.M.Stuart}, {\em {A} {N}ote on {D}iffusion {L}imits of
  {C}haotic {S}kew {P}roduct {F}lows}, Nonlinearity, 24 (2011), pp.~1361--1367.

\bibitem{Niederreiter1992}
{\sc H.~Niederreiter}, {\em {R}andom {N}umber {G}eneration and {Q}uasi-{M}onte
  {C}arlo {M}ethods}, CBMS-NSF regional conference series in applied
  mathematics, Society for Industrial and Applied Mathematics, 1992.

\bibitem{Nolen2012}
{\sc J.~Nolen, G.~A. Pavliotis, and A.~M. Stuart}, {\em {M}ultiscale
  {M}odelling and {I}nverse {P}roblems}, in Numerical Analysis of Multiscale
  Problems, O.~L. I.G.~Graham, T.Y.~Hou and R.~Scheichl, eds., vol.~83,
  Springer, 2012, pp.~1--34.

\bibitem{Papavasiliou2009}
{\sc A.~Papavasiliou, G.~A. Pavliotis, and A.~M. Stuart}, {\em Maximum
  likelihood drift estimation for multiscale diffusions}, Stochastic Process.
  Appl., 119 (2009), pp.~3173--3210.

\bibitem{Pavliotis2007}
{\sc G.~A. Pavliotis and A.~M. Stuart}, {\em {P}arameter {E}stimation for
  {M}ultiscale {D}iffusions}, J. Stat. Phys., 127 (2007), pp.~741--781.

\bibitem{Pavliotis2008book}
\leavevmode\vrule height 2pt depth -1.6pt width 23pt, {\em {M}ultiscale
  {M}ethods: {A}veraging and {H}omogenization}, Springer, 2008.

\bibitem{Penrose1956}
{\sc R.~Penrose}, {\em On best approximate solutions of linear matrix
  equations}, Math. Proc. Cambridge Philos. Soc., 52 (1956), pp.~17--19.

\bibitem{Pradas2012}
{\sc M.~Pradas, G.~A. Pavliotis, S.~Kalliadasis, D.~T. Papageorgiou, and
  D.~Tseluiko}, {\em Additive noise effects in active nonlinear spatially
  extended systems}, European J. Appl. Math., 23 (2012), pp.~563--591.

\bibitem{Pradas2011}
{\sc M.~Pradas, D.~Tseluiko, S.~Kalliadasis, D.~T. Papageorgiou, and G.~A.
  Pavliotis}, {\em {N}oise {I}nduced {S}tate {T}ransitions, {I}ntermittency,
  and {U}niversality in the {N}oisy {K}uramoto-{S}ivashinksy {E}quation}, Phys.
  Rev. Lett., 106 (2011), pp.~060602--1--060602--4.

\bibitem{PrakasaRao1999}
{\sc B.~L.~S. {Prakasa Rao}}, {\em {S}tatistical {I}nference for {D}iffusion
  {T}ype {P}rocesses}, vol.~8 of Kendall's Library of Statistics, Arnold, 1999.

\bibitem{Savva2010}
{\sc N.~Savva, S.~Kalliadasis, and G.~A. Pavliotis}, {\em {T}wo-{D}imensional
  {D}roplet {S}preading over {R}andom {T}opographical {S}ubstrates}, Phys. Rev.
  Lett., 104 (2010), pp.~084501--1--084501--4.

\bibitem{Theodoropoulos2000}
{\sc C.~Theodoropoulos, Y.-H. Qian, and I.~G. Kevrekidis}, {\em {"Coarse"
  Stability and Bifurcation Analysis Using Time-Steppers: A Reaction-Diffusion
  Example}}, Proc. Natl. Acad. Sci. USA, 97 (2000), pp.~9840--9843.

\bibitem{Vanden-Eijnden2003}
{\sc E.~Vanden-Eijnden}, {\em Numerical techniques for multi-scale dynamical
  systems with stochastic effects}, Commun. Math. Sci., 1 (2003), pp.~385--391.

\bibitem{Zhang2005}
{\sc L.~Zhang, P.~A. Mykland, and Y.~A{\"{\i}}t-Sahalia}, {\em A {T}ale of
  {T}wo {T}ime {S}cales: {D}etermining {I}ntegrated {V}olatility with {N}oisy
  {H}igh-{F}requency {D}ata}, J. Amer. Statist. Assoc., 100 (2005),
  pp.~1394--1411.

\bibitem{Zwanzig2001}
{\sc R.~Zwanzig}, {\em {N}onequilibrium {S}tatistical {M}echanics}, Oxford
  University Press, 2001.

\end{thebibliography}
